\documentclass[12pt,A4]{article}
\usepackage{amssymb, amsmath, enumerate, geometry,epsf}
\usepackage[all]{xy}
\begin{document}

\newtheorem{thm}{Theorem}
\newtheorem{lm}[thm]{Lemma}
\newtheorem{prop}[thm]{Proposition}
\newtheorem{cor}[thm]{Corollary}
\newtheorem{defi}[thm]{Definition}

\numberwithin{equation}{section}

\newcommand{\proof}{{\bf Proof\ }}
\newcommand{\eproof}{$\blacksquare$\bigskip}
\newcommand{\CC}{\mathcal C}
\newcommand{\CF}{\mathcal F}
\newcommand{\CR}{\mathcal R}
\newcommand{\C}{\Bbb C}
\newcommand{\tens}{\otimes}
\renewcommand{\o}{{}_{\scriptscriptstyle(1)}}
\renewcommand{\t}{{}_{\scriptscriptstyle(2)}}
\newcommand{\thr}{{}_{\scriptscriptstyle(3)}}
\newcommand{\fo}{{}_{\scriptscriptstyle(4)}}
\newcommand{\Nat}{{\mathrm{Nat}}}
\newcommand{\id}{{\mathrm{id}}}
\newcommand{\isom}{{\cong}}
 
\newcommand{\prob}{\mathop{\mathrm{Pr}}}
\newcommand{\eps}{\varepsilon}
\newcommand{\raction}{\triangleleft}
\newcommand{\Del}{\Delta}
\newcommand{\del}{\delta}
\newcommand{\Lan}{\Lambda}
\newcommand{\Gam}{\Gamma}
\newcommand{\gam}{\gamma}
\newcommand{\ten}{\bigotimes}
\newcommand{\sten}{\otimes}
\newcommand{\Om}{\Omega}
\newcommand{\sfi}{\varphi}
\newcommand{\om}{\omega}
\newcommand{\lan}{\lambda}
\newcommand{\sig}{\sigma}
\newcommand{\al}{\alpha}
\newcommand{\all}{\forall}
\newcommand{\flsh}{\longrightarrow}
\newcommand{\ld}{\ldots}
\newcommand{\cd}{\cdots}
\newcommand{\bu}{\bullet}
\newcommand{\eqn}[2]{\begin{equation}#2\label{#1}\end{equation}}

\textheight 23.6cm
\textwidth 16cm
\topmargin -.2in \headheight 0in \headsep 0in
\oddsidemargin 0in \evensidemargin 0in
\topskip 28pt

\title{Braided Cyclic Cocycles and Non-Associative Geometry}

\author{S.E. Akrami\footnote{Supported financialy by British Council and Ministry of Science of Iran}\\Department of Mathematics, University of Tehran, Iran\\{+}\\
S. Majid \\ School of Mathematics, Queen Mary, University of London\\
London E14NS, UK}

\date{31 May, 2004}

\maketitle

\begin{abstract}
We use monoidal category methods to study the noncommutative geometry of nonassociative algebras obtained by a Drinfeld-type cochain twist. These are the so-called quasialgebras and include the octonions as braided-commutative but nonassociative coordinate rings, as well as quasialgebra versions $\CC_{q}(G)$ of the standard $q$-deformation quantum groups. We introduce the notion of ribbon algebras in the category, which are algebras equipped with a suitable generalised automorphism $\sigma$, and obtain the required generalisation of cyclic cohomology. We show that this \emph{braided cyclic cocohomology} is invariant under a cochain twist. We also extend to our generalisation the relation between cyclic cohomology and differential calculus on the ribbon quasialgebra. The paper includes differential calculus and cyclic cocycles on the octonions as a finite nonassociative geometry, as well as the algebraic noncommutative torus as an associative example. 
\end{abstract}
\section{Introduction}

In the influential work  \cite{Drinfeld}, Drinfeld extended quantum groups or Hopf algebras to a more general notion of quasi-Hopf algebras stable under conjugation of the coproduct by a 'twist'. In a dual form the cotwist element $F$ is a cochain and modifies the product of the coquasiHopf algebra. In \cite{Ma:tan} this construction was formulated as a  monoidal equivalence between the comodule category of the coquasiHopf algebra $H$ and that of the cotwisted coquasiHopf algebra $H^{F}$. 
The differential geometry of quantum groups $H^{F}$ from this point of view and assuming $F$ was a cocycle (so that we stay in the associative setting) was studied in \cite{MaOeck}. The more general coquasiHopf setting was used recently in \cite{BegMa} and applied to the standard quantum groups $\C_{q}(G)$. This work proved that there is no associative differential algebra on the standard quantum groups with classical dimensions (i.e. deforming the classical case in a strict sense) but that this is possible as a supercoquasi-Hopf algebra $\Omega(\C_{q}(G))$. This could be considered a first hint that nonassociative geometry is necessary for a full understanding even of ordinary quantum groups. It also suggests that  one should take seriously nonassociative coordinate algebras themselves (not just their exterior algebras) and moreover in much greater generality than just (coquasi)Hopf algebras alone.  

We do this in the present paper for algebras $A$ in monoidal Abelian categories. The idea is that the nicest nonassociative algebras, which we call {\em quasialgebras}, should be ones which are nonassociative but which may be viewed as associative by deforming the tensor product to a monoidal category with nontrivial associator $\Phi_{{U,V,W}}:U\tens (V\tens W)\to (U\tens V)\tens W$ for the rebracketting of tensor products of objects $U,V,W$. This complements the established idea of using braided categories to view certain noncommutative algebras as 'commutative' with respect to a nontrivial braiding $\Psi_{{V,W}}:V\tens W\to W\tens V$ for objects $V,W$, see \cite{Ma:book}. Similarly, building on work of Drinfeld\cite{Drinfeld} we will find quasialgebra versions $\CC_{q}(G)$ of the standard quantum groups which are nonassociative but more commutative (i.e. one can trade one feature for the other). 

In general terms, we study several key constructions borrowed from noncommutative geometry\cite{Con}\cite{KMT} but now in the setting of quasialgebras, i.e. of algebras in monoidal categories with nontrivial associator.  The idea is to think of the quasialgebra geometrically as by definition the coordinates of a 'quasiassociative' space, which may also be noncommutative (for example, commutative with respect to a nontrivial braiding) and hence a 'quantum quasispace'. Thus, in Section~3 we will associate a cocyclic module (see\cite{Con}) to any algebra in a braided monoidal Ab-category. We show that the morphism 
\[ \lan=(\sig\sten id)\Psi\]  provides for us a cyclicity morphism in this category and hence a 'braided cyclic cohomology' theory. Here $\Psi:A^{\sten (n+1)}\sten A\flsh A\sten A^{\sten(n+1)}$ is the braiding isomorphism and $\sig:A\flsh A$ is a ``generalised algebra automorphism''. For reasons which will become clear,  we call it a {\em ribbon automophism} on the algebra $A$ and call $(A,\sigma)$ a {\em ribbon (quasi)algebra}. In the associative trivially braided case $\sigma$ becomes an automorphism in the usual sense and our braided cyclic cohomology reduces to the 'twisted' cohomology in \cite{KMT} which has been successfully applied to quantum groups such as $\C_{q}[SL_{2}]$. We also study how  braided cyclic cocycles relate to differential calculus in the monoidal category.

The other key goal of the paper concerns the provision of examples. Indeed, the need for some kind of nonassociative geometry is hinted at from several directions in mathematical physics including string theory. Its need is also indicated from Poisson geometry, where the idea of a generalised Poisson bracket violating the usual Jacobi identity is proposed\cite{AleKos}. It turns out that an adequate class that appears to cover such examples is based on the use of Drinfeld type cotwists, but not for (coquasi)Hopf algebras $H$ as above. Rather,  we consider an algebra $A$ in the  monoidal category of $H$-comodules. After applying the monoidal equivalence one has an algebra $A_{F}$ in the category of $H^{F}$- comodules.  Indeed, all algebras and algebraic constructions 'gauge transform' in this way. The $A_{F}$ construction was introduced in \cite{Ma:euc} (in a module version) and for our purposes  takes the form\cite{Ma:book}
\[ A_{F}=A,\quad {\rm with\ the\ new\ product}\quad a._{F} b=F(a\o,b\o)a\t b\t\]
where $a\o\tens a\t$ denotes the (left) coaction. It turns out that a great many noncommutative algebras  of interest fit into this cotwist framework for suitable $H$ and $F$. Indeed, switching on $F$ is a useful formulation of quantisation as an extension of the Moyal product: 
Even if $A$ is commutative, $A_{F}$ in general becomes noncommutative when $F$ is not symmetric since
\[a._{F} b=F(a_{(1)},b_{(1)})F^{-1}(b_{(2)},a_{(2)})b_{(3)}._{F} a_{(3)}. \]
More relevant for us, even if $A$ is associative, $A_{F}$ in general becomes nonassociative unless $F$ obeys a 2-cocycle condition\cite{Ma:book}. Recent applications include \cite{Sit}\cite{Var}\cite{Ma:cli} in the associative case and \cite{AlbMa}, in the nonassociative case.  

 Section~4 contains theorems about how  braided cyclic cohomology and differential geometry respond under such cotwists.  Thus, suppose for the sake of discussion that $A$ possesses a left covariant differential calculus, $\Om=\bigoplus_{k=0}^n\Om_k$ which is super commutative i.e. two homogeneous differential forms $\om$ and $\om'$ commute up to a sign $(-1)^{|\om||\om'|}$. Thus in particular functions (0-forms) and $n$-forms commute. Therefore if $\int:\Om_n\flsh\mathbb{C}$ is a closed graded trace in the sense of Connes  then its character is a cyclic cocycle\cite{Con}. Now if we cotwist the superalgebra of differential forms with the same cochain from one side as above, we obtain a calculus $\Om_F$ for $A_{F}$ but now functions and $n$-forms no longer commute. The  noncommutativity is controlled by
$$\om\om'=(-1)^{|\om||\om'|}F(\om_{(1)},\om'_{(1)})F^{-1}(\om'_{(2)},\om_{(2)})\om'_{(3)} \om_{(3)}$$
where the product is the cotwisted $._{F}$ one. Consequently the character $\sfi(a^0,\ldots,a^n)=\int a^0da^1\cd da^n$ after cotwisting is no longer a cyclic cocycle but obeys\begin{eqnarray*}
&&\sum_{i=0}^n(-1)^i\sfi (a^0,\ldots,a^i a^{i+1},\ldots,a^{n+1})+\\
&&\quad +(-1)^{n+1}F(a^0_{(1)}\cd a^n_{(1)},a^{n+1}_{(1)})F^{-1}(a^{n+1}_{(2)},a^0_{(2)}\cd a^n_{(2)})\sfi(a^{n+1}_{(3)}a^0_{(3)},a^1_{(3)},\ld, a^n_{(3)})=0. \end{eqnarray*}
The corresponding formula in the case where $F$ is just a cochain not a cocycle, is much more involved and contains the associator in its formula, (see Section~4.1). We obtain, rather, a braided cyclic cocycle in the 'gauge equivalent' category of $H^{F}$-comodules. Section~4.1 also contains rather concrete formulae when the background Hopf algebra $H$ is coquasitriangular. Section~4.2 specialises the theory to the important case where  in fact $H$ is the group algebra of an abelian group, which is the setting needed for many examples including the octonions.
  
 Finally, Section~5 presents a collection of key examples demonstrating the theory of paper. We explicitly give the differential calculus and a cyclic cocycle on the octonions as a finite nonassociative geometry, as well as the usual (algebraic) noncommutative torus. Section~5.3 also outlines the theory applied to formal deformation theory,  where we obtain the quasialgebras $\CC_{q}(G)$ as mentioned above, using Drinfeld's associator obtained from solving the Knizhnik-Zamalochikov equations. We can also in principle use cotwising to deformation-quantise  the quasiPoisson manifold structure on a $G$-manifold $M$ proposed in \cite{AleKos}  (which was not achieved before), which we do as a quasialgebra $\CC_{q}[M]$. We also construct the differential calculus and braided cyclic cocycles on all these quasialgebras.  Further details of these potential examples will be presented elsewhere.
     
On the technical side, we start in the preliminary Section 2, by explicitly embedding, in a canonical way, any general (relaxed) braided monoidal Ab-category in a strict braided monoidal Ab-category. This underlies Mac Lane's coherence theorem\cite{Mac} and ensures that one can work with a relaxed category like a strict one. We also recall the Drinfeld's ``gauge transformation'' (in a dual cotwist sense) at the level of braided monoidal Ab-categories as in \cite{Ma:tan}. 

We conclude the introduction  with the geometric motivation behind our theory. In what follows we will consider `branched ribbon tangles', a modification of the usual notion of ribbon tangles (see \cite{Turaev} and the references there). We are not going to give a precise meaning of a branched ribbon tangle, but limit ourselves to an informal discussion. Thus let us define $d_i,~0\le i\le n $ and $d_{n+1}$ to be the isotopy type of the  branched ribbon tangles in the strip $\mathbb{R}^2\times [0,1]$ in Figure~1, (a) and (b) respectively. Here by isotopy we mean isotopy in $\mathbb{R}^2\times[0,1]$ constant in boundary intervals on the lines $z=0$ and $z=1$ in the plane $x=0$. Then intuitively we have the isotopies in Figure~1 parts (c), (d), (e) and (f). In part (d) we used the isotopy in part (g). All these isotopies should be clear except (f) which may need more explanation; in fact ignoring vertical bands indexed from 1 to n+1 in part (f), the left hand side of (f) is isotopic with the left hand side of (h). Now considering the last diagram in (h), rotate the upper branch in this diagram by 360 degree to give the right hand side of (f).
Now if, as in \cite{Turaev} or \cite{Ma:book} we define the composition of two branched ribbon tangles by putting one on  top of the other one and compressing the resulting diagram to the strand $\mathbb{R}^2\times[0,1]$. Then  parts  (c), (d), (e) and (f) of Figure~1 mean
$$d_id_j=d_{j-1}d_i$$
\begin{figure}
\[ \epsfbox{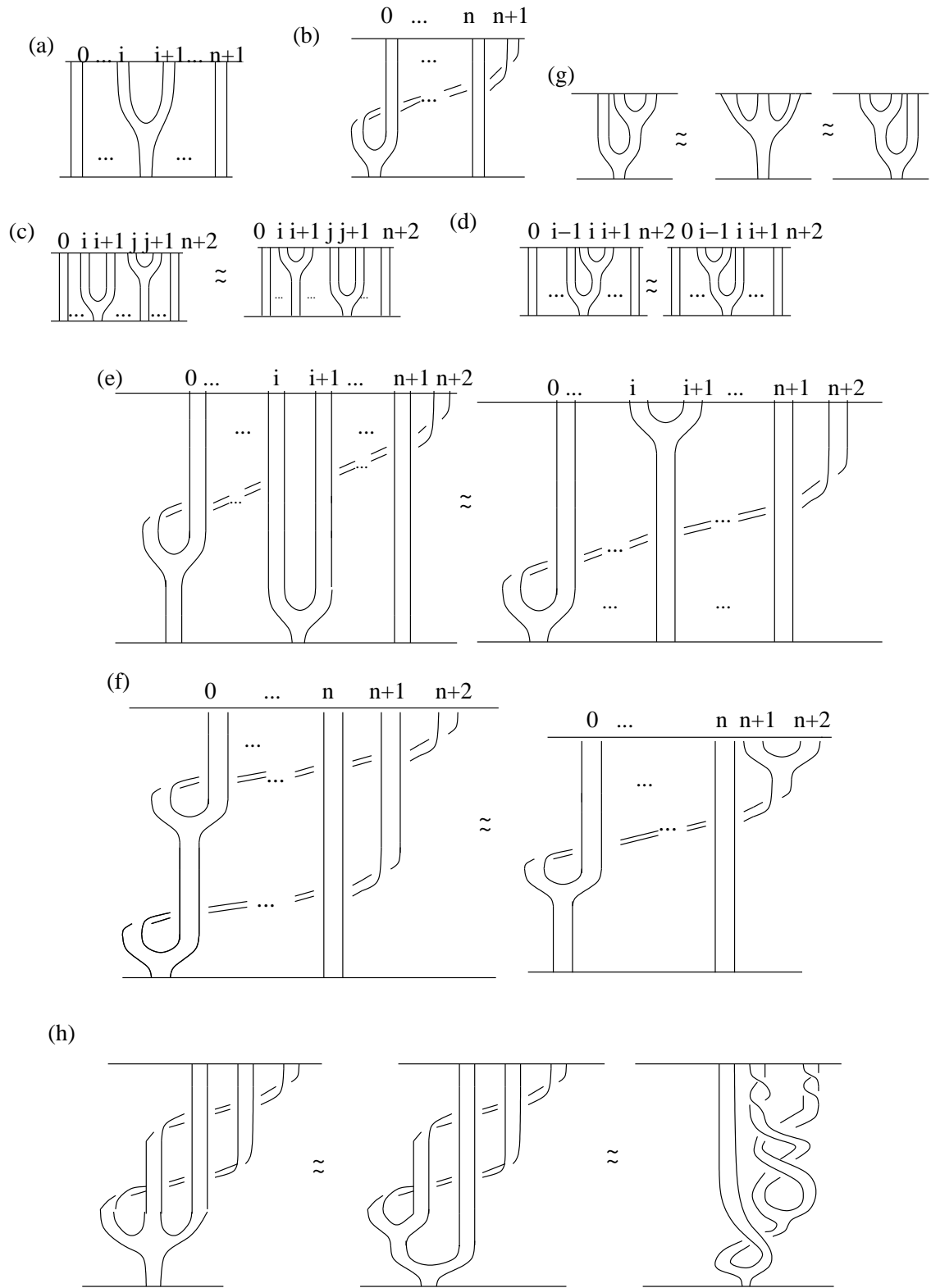} \]
\caption{(a),$ 0\le i\le n$. (c), $0\le i<j-1\le n$. (d), $1\le j \le n+1$. (e),$ 0\le i\le n$.  }
\end{figure}

These relations are the main content of the notion of a simplicial object in the category. Note that one  can work essentially in the category of ribbon tangles as in\cite{Turaev}, where roughly speaking, a morphism is just an isotopy type of a ribbon tangle. This would need, however, some geometric considerations whereas we prefer to work in a purely algebraic manner. Thus in Section~3 we  axiomatize the precise assumptions which will lead us to above relations in the context of a  general braided monoidal category. In this case, for braided categories we use the usual graphical calculus, which should not be confused with the above isotopy argument on actual ribbon graphs which are embedded surfaces in space. We recall that any braided category is the image of the category of braids which allows for the representation of the axiomatic properties of  $\tens,\Psi$ by strands and their braiding in the graphical notation (see the Preliminaries section).  In the same way, one can define a functor from the category of branched ribbon graphs, (see\cite{Turaev}) to our category in Section~3 such that two isotopic branched ribbon tangles have the same image under this functor, and the image of branched ribbon tangles under this functor can be considered as graphical symbols.  

\textbf{Acknowledgment}. This paper was written during the year-long visit by the first author to Queen Mary. He would like to thank the staff and students there for their kindness  and the British Council and Iranian Ministry of Science for making the visit possible.

\section{Preliminaries: monoidal categories, coherence and `gauge' equivalence}

Here we establish the basic notations and methods needed in the paper. In order to be able to work effectively with the nonassociativity, we need in particular to fix some conventions on bracketting and explain related issues. Then we recall the notion of monoidal equivalence used in the construction of our examples. 

 Let  $(\mathcal{C},\sten,\Phi,\Psi )$ be a braided monoidal Ab-category, where $\Phi $ is the associator and $\Psi$ is the braiding. We refer to \cite{Turaev,Ma:book} for the axioms, but briefly $\tens:\CC\times\CC\to \CC$ is a functor and $\Phi\in\Nat(\ \tens\ (\ \tens\ ), (\ \tens\ )\tens\ )$ and $\Psi\in\Nat(\tens,\tens^{{\rm op}})$  are natural isomorphisms obeying respectively the pentagon and two hexagon identities with regard to composite objects. The pentagon identifies two ways to go from $U\tens(V\tens(W\tens Z))\to ((U\tens V)\tens W)\tens Z$ using two and three applications of $\Phi$ respectively. The hexagons identify $(V\tens W)\tens Z\to Z\tens (V\tens W)$ in two ways and $V\tens (W\tens Z)\to (W\tens Z)\tens V$ in two ways, using $\Psi,\Phi$. The coherence theorem of Mac Lane in the symmetric case (when $\Psi^{2}=\id$ is assumed) and of Joyal and Street\cite{JoyStr} in the general braided case, ensures that once these are assumed, one may (a) drop brackets entirely, inserting $\Phi$ as needed for any composition to make sense, and (b) represent $\Psi,\Psi^{{-1}}$ by braid crossings $\epsfbox{braid.eps}$ and inverse braid crossings $\epsfbox{braidinv.eps}$; if two braid compositions are equal then so are the compositions of $\Psi,\Psi^{{-1}}$. We assume the monoidal category is unital and denote the unit object by \mbox{\boldmath $1$}, where for simplicity we will assume that $U\sten \mbox{\boldmath $1$}= \mbox{\boldmath $1$}\sten U=U,~ \all U\in\mathcal{C}$ and $\Phi_{U,V,W}=id$ whenever at least one of the objects $U,V$ or $W$ is \mbox{\boldmath $1$}.   Similarly for the braiding when present. By Ab-monoidal category we mean that for each pair of objects $U$ and $V$ the set $Hom (U,V)$ is an additive abelian group such that the composition and tensor product of morphisms are bilinear. Then one can easily show that $K:=Hom(\mbox{\boldmath $1$},\mbox{\boldmath $1$})$ is a commutative ring with unit.

Next, in order to be completely explicit, we canonically extend  any  braided monoidal Ab-category $\mathcal{C}$ into a strict braided monoidal Ab-category.  First of all, using induction on $n$ we define a family $\Lan_n$ of sets  as follows.  Let for $n>0$, $\mathcal{C}^n$ be the $n$ times Cartesian product of  $\mathcal{C}$ with itself. Let   $\Lan_1$ be the single set whose element is the identity functor on $\mathcal{C}$ and let $\Lan_2$ be the single set whose element is the functor $\sten:\mathcal{C}\times\mathcal{C}\flsh\mathcal{C}$. Now suppose that  $\Lan_k, k< n$ has been defined already such that for each $\al\in\Lan_k$ there is a associated functor $\bar{\al}:\mathcal{C}^k\flsh\mathcal{C}$, then we define  $\Lan_n$ for $n>2$ to be the set of all pairs $(\al,\beta),~\al\in\Lan_k,~\beta\in\Lan_l,~k<n,l<n,~k+l=n$ and we associate to the pair $(\al,\beta)$ the functor $\sten(\al\times\beta) :\mathcal{C}^n\flsh\mathcal{C}$. For  $ U_i\in\mathcal{C},1\le i\le n$ and $\al\in\Lan_n$ the symbol $(U_1\sten\cd\sten U_n;\al)$ or $U_1\sten_\al\cd\sten_\al U_n$ will denote the object $\bar{\al}(U_1,\ld,U_n)$ .

Now let $\bar{\mathcal{C}}$ be the category whose objects are  symbols  $U_1\sten\cd\sten U_n,~n=1,2,\ld, U_i\in \mathcal{C}$ such that if $n>1$ then $U_i\ne \mbox{\boldmath $1$}, \all i$. The object $ \mbox{\boldmath $1$}$ of $\bar{\mathcal{C}}$ is called the degree zero object and the objects $U_1\sten\cd\sten U_n,~n=1,2,\ld, \mbox{\boldmath $1$}\ne U_i\in \mathcal{C}$ are called of degree $n$.  To define morphisms in this category we recall that by Mac Lane's coherence theorem there exists a unique family of isomorphisms $I^\al_\beta=I^\al_\beta(U_1,\ld, U_m):(U_1\sten\cd\sten U_m;\al)\flsh(U_1\sten\cd\sten U_m;\beta),~\al,\beta\in\Lan_m,\mbox{\boldmath $1$}\ne U_i\in\mathcal{C}$, induced by the associator satisfying  $I^\al_\beta=id$ for $m=1,2$ and
\eqn{}{I^\beta_\gam I^\al_\beta=I^\al_\gam,~~~\all\al,\beta,\gam\in\Lan_m}
\eqn{}{I^{\al_1}_{\al_2}(U_1,\ld,U_k)\sten I^{\beta_1}_{\beta_2}(U_{k+1},\ld,U_m)=I^{(\al_1,\beta_1)}_{(\al_2,\beta_2)}(U_1,\ld,U_m),~~~\al_1,\al_2\in\Lan_k,\beta_1,\beta_2\in\Lan_{m-k}}
\eqn{}{I^{(id,\sten)}_{(\sten,id)}(U_1,U_2,U_3)=\Phi_{U_1,U_2,U_3}}
We call these isomorphisms , the \emph{higher degree associators} of  $\mathcal{C}$.

Next we define an equivalence relation among the morphisms in $\bigcup_{\al\in \Lan_m,\beta\in \Lan_n}Hom_\mathcal{C}((U_1\sten\cd\sten U_m;\al),(V_1\sten\cd\sten V_n;\beta))$ by
 \begin{eqnarray}
f_1\sim f_2 \Longleftrightarrow f_2=I^{\beta_1}_{\beta_2}(V_1,\ld,V_n) f_1 I^{\al_2}_{\al_1}(U_1,\ld,U_m)
\end{eqnarray}
for $f_i:(U_1\sten\cd\sten U_m;\al_i)\flsh(V_1\sten\cd\sten V_n;\beta_i)), ~\al_i\in \Lan_m, \beta_i\in \Lan_n,i=1,2$. One can easily show that the above relation is an equivalence relation and each equivalence class contains one and only one representative in each set $Hom_\mathcal{C}((U_1\sten\cd\sten U_m;\al),(V_1\sten\cd\sten V_n;\beta))$. We now define $Hom_{\bar{\mathcal{C}}}(U_1\sten\cd\sten U_m,V_1\sten\cd\sten V_n)$ be the equivalence classes of all morphisms in the set $\bigcup_{\al\in \Lan_m,\beta\in \Lan_n}Hom_\mathcal{C}((U_1\sten\cd\sten U_m;\al),(V_1\sten\cd\sten V_n;\beta))$. We denote the equivalence classes by notation $[f]$.

For morphisms $f: U_1\sten\cdots\sten U_l\flsh V_1\sten\cdots\sten V_m$ and $g:V_1\sten\cdots\sten V_m\flsh W_1\sten\cdots\sten W_n$  in $\bar{\mathcal{C}}$ we define composition $gf$ to be the class of morphism $g^\beta_\gam f^\al_\beta$. Where $ f^\al_\beta$ and $g^\beta_\gam$ are the representatives of $f$ and $g$ in $Hom_\mathcal{C}((U_1\sten\cd\sten U_l;\al),(V_1\sten\cd\sten V_m;\beta))$ and $Hom_\mathcal{C}((V_1\sten\cd\sten V_m;\beta),(W_1\sten\cd\sten W_n;\gam))$ respectively, for arbitrary $\al\in\Lan_l,\beta\in\Lan_m,\gam\in\Lan_n$. One can easily show that this composition is well-defined and is associative. And the classes of identity morphisms in $\mathcal{C}$ are identity morphisms in $\bar{\mathcal{C}}$.

To define a monoidal structure on $\bar{\mathcal{C}}$ we define $\mbox{\boldmath $1$}\sten\mbox{\boldmath $1$}:=\mbox{\boldmath $1$}$ and $\mbox{\boldmath $1$}\sten (U_1\sten\cdots\sten U_m):=U_1\sten\cdots\sten U_m=:(U_1\sten\cdots\sten U_m)\sten\mbox{\boldmath $1$}$ and $(U_1\sten\cdots\sten U_m)\sten(V_1\sten\cdots\sten V_n):=U_1\sten\cdots\sten U_m\sten V_1\sten\cdots\sten V_n$ for $\mbox{\boldmath $1$}\ne U_i,V_i$ . And  for $f:U_1\sten\cdots\sten U_m\flsh V_1\sten\cdots\sten V_n $ and $f': U'_1\sten\cdots\sten U'_{m'}\flsh  V'_1\sten\cdots\sten V'_{n'}$ we define $f\sten f'$ to be the equivalence class of the morphism $f^\al_\beta\sten {f'}^{\al'}_{\beta'}:(U_1\sten\cdots\sten U_m\sten U'_1\sten\cdots\sten U'_{m'};(\al,\al'))\flsh (V_1\sten\cdots\sten V_n\sten V'_1\sten\cdots\sten V'_{n'};(\beta,\beta')),~\al\in\Lan_m,\al'\in\Lan_{m'},\beta\in\Lan_n,\beta\in\Lan_{n'}$. Again one can show this is well-defined by using (2.2). From naturality of the  associator, it is obvious that this tensor product is associative and the unit object of $\mathcal{C}$ is also unit object of this product.
Therefore $\bar{\mathcal{C}}$ is a strict monoidal category. Similarly the addition of two morphisms $f,g:U_1\sten\cdots\sten U_m\flsh V_1\sten\cdots\sten V_n $ is defined to be the equivalence class of the morphism $f^\al_\beta+g^\al_\beta,~\al\in\Lan_m,~\beta\in\Lan_n$, and it is straightforward to show this is well-defined and   $\bar{\mathcal{C}}$ is a strict monoidal Ab-category.

Finally, when $\CC$ is braided, we define braid isomorphisms $\bar{\Psi}_{U_1\sten\cdots\sten U_m,V_1\sten\cdots\sten V_n}$ to be the equivalence class of morphism $\Psi_{(U_1\sten\cdots\sten U_m;\al),(V_1\sten\cdots\sten V_n;\beta)}$. It is easy to show that $\bar{\mathcal{C}}$ equipped with this isomorphisms becomes a strict braided monoidal Ab-category.

Note that the strict braided monoidal category  $\bar{\mathcal{C}}$ is an extension of  $\mathcal{C}$ as an Ab-category but not as a braided monoidal one. Nevertheless, the braided monoidal structure of $\CC$ is used in building $\bar\CC$ and allows us to replace any non-strict $\CC$ by the strict  $\bar\CC$ by regarding its objects and morphisms in $\bar\CC$. Note that  $\mathcal{C}$ does not inherit an associative tensor product from $\bar{\mathcal{C}}$ because if we regard two objects $U$ and $V$ in  $\mathcal{C}$ as first degree objects of $\bar{\mathcal{C}}$ then their tensor product $U\sten V$ in $\bar{\mathcal{C}}$ is a second degree object of $\bar{\mathcal{C}}$ and therefore it is not in $\mathcal{C}$ which is the set of objects of degree less than two in $\bar{\mathcal{C}}$ ). For example we know that an algebra in a strict monoidal category is defined to be an object $A$ equipped with a morphism $m:A\sten A\flsh A$ such that $m(id\sten m)=m(m\sten id)$ as a morphism from $A\sten A\sten A$ to $A$. Now let $\mathcal{C}$ be a non-strict monoidal category, but view it in $\bar\CC$.    We define similarly, an algebra $A$ to be an object in $\mathcal{C}$ ( a first degree object in  $\bar{\mathcal{C}}$) equipped with a morphism   $m:A\sten A\flsh A$ ( a morphism from a second degree object to a first degree object in $\bar{\mathcal{C}}$) such that $m(id\sten m)=m(m\sten id)$ as an equality between morphisms from $A\sten A\sten A$ ( a third degree object in $\bar{\mathcal{C}}$) to $A$ ( a first degree object in $\bar{\mathcal{C}}$). If in addition there exists a morphism $\eta:\mbox{\boldmath $1$}\flsh A$ such that $m(id\sten\eta)=m(\eta\sten id)=id$ then we call $(A,m,\eta)$ a unital algebra in  $\mathcal{C}$.

The above construction of the family $\{\Lan_n\}_n$ can be done for any set  equipped with a binary action on it. Explicitly let $A$ be a set  with  a binary action $.:A \times A \flsh A$, similarly, using induction, for each $n$ we have a set of $n$-fold  actions $A\times \cd\times A\flsh A$ which we denote this set again by $\{\Lan_n\}_n$ and for $\al\in \Lan_n$ and $a^i\in A,~1\le i\le n$ we use the notation $a^1._\al\cd._\al a^n$ for $\al(a^1,\ld,a^n)$.\\ 

We also want to recall quite explicitly the notion of a tensor functor between braided monoidal categories. It is known\cite{Ma:tan} that equivalence by such tensor functors is the correct notion of 'gauge transformation' relevant to the Drinfeld cotwist. At the moment, we give the general categorical setting for this. Thus, let $(\mathcal{C},\sten,\Phi,\Psi)$ and $(\mathcal{C}',\sten',\Phi',\Psi')$ be two braided unital monoidal Ab-categories  and let $T:\mathcal{C}\flsh \mathcal{C}'$ be an additive functor such that $T(\mbox{\boldmath $1$})=\mbox{\boldmath $1$}'$ and suppose that there exists a natural isomorphism $\mathcal{F}$ between the functors $\sten'(T\times T),T\sten :\mathcal{C}\times \mathcal{C}\flsh \mathcal{C}' $, i.e. there exist a family of isomorphisms $\mathcal{F}_{U,V}:T(U)\sten'T(V)\flsh T(U\sten V)$ in $\mathcal{C}'$, $\all U,V\in\mathcal{C}$, such that $\mathcal{F}_{U,\mbox{\boldmath $1$}}=\mathcal{F}_{\mbox{\boldmath $1$},U}=id_{T(U)}$ and for all objects $U_i,V_i, i=1,2$ and morphisms $f:U_1\flsh U_2,~g:V_1\flsh V_2 $ in $\mathcal{C}$ we have 
\eqn{tf}{T(f\sten g)\mathcal{F}_{U_1,V_1}=\mathcal{F}_{U_2,V_2}(T(f)\sten'T(g))}
Now suppose that    
\eqn{tf-phi}{T(\Phi_{U,V,W})\mathcal{F}_{U,V\sten W}(id_{T(U)}\sten'\mathcal{F}_{V,W})=\mathcal{F}_{U\sten V,W}(\mathcal{F}_{U,V}\sten id_{T(W)})\Phi'_{T(U),T(V),T(W)}}
and
\eqn{tf-psi}{T(\Psi_{U,V})\mathcal{F}_{U,V}=\mathcal{F}_{V,U}\Psi'_{T(U),T(V)}}
$\all U,V,W \in \mathcal{C}$.

\begin{defi}\label{tenfun} A tensor functor   between two braided monoidal Ab-categories  $(\mathcal{C},\sten,\Phi,\Psi)$ and $(\mathcal{C}',\sten',\Phi',\Psi')$ is a pair $(T,\mathcal{F})$ as above. A 'gauge transformation' between  braided monoidal Ab-categories is an invertible tensor functor, in which case we say that the categories are 'gauge equivalent'.  \end{defi}

 Indeed, if  $T$ is an invertible functor, we set  
$${\mathcal{F}'}_{U',V'}:=T^{-1}(\mathcal{F}^{-1}_{T^{-1}(U'),T^{-1}(V')}),~~~\all U',V'\in \mathcal{C}'$$ 
Then $(T^{-1},{\mathcal{F}'})$ is a tensor functor from  $ (\mathcal{C}',\sten',\Phi',\Psi')$ to $(\mathcal{C},\sten,\Phi,\Psi)$. To prove this let $U'_i,V'_i,U',V'$ and $W'$  be objects in $ \mathcal{C}'$, where $ i=1,2$, and let $f':U'_1\flsh U'_2,~g':V'_1\flsh V'_2$ be morphisms in $ \mathcal{C}'$, then we set $U_i:=T^{-1}(U'_i),V_i:=T^{-1}(V'_i),~i=1,2,~f:=T^{-1}(f')$ and $g:=T^{-1}(g')$. Then from (2.5) we get $\mathcal{F}^{-1}_{U_2,V_2}T(f\sten g)=(f'\sten'g')\mathcal{F}^{-1}_{U_1,V_1}$. Thus applying $ T^{-1}$ to this relation we get
\eqn{}{T^{-1}(f'\sten' g'){\mathcal{F}'}_{U'_1,V'_1}={\mathcal{F}'}_{U'_2,V'_2}(T^{-1}(f')\sten'T^{-1}(g'))}
which is counterpart of (2.5) for pair $(T^{-1},{\mathcal{F}'})$. The counterpart of (2.7) for the pair $(T^{-1},{\mathcal{F}'})$ can be proved similarly. Let us prove (2.6) for $(T^{-1},{\mathcal{F}'})$. At first note that if we apply (2.8) for $f'=id_{U'},g'=\mathcal{F}^{-1}_{V,W}$ then we get $T^{-1}(id_{U'}\sten'\mathcal{F}^{-1}_{V,W} ){\mathcal{F}'}_{U',T(V\sten W)}={\mathcal{F}'}_{U',{V'\sten W'}}(id_{T^{-1}(U')}\sten'T^{-1}(\mathcal{F}^{-1}_{V,W}))={\mathcal{F}'}_{U',{V'\sten W'}}(id_{T^{-1}(U')}\sten'{\mathcal{F}'}_{V',W'})$. And similarly applying (2.8) for $f'=\mathcal{F}^{-1}_{U,V}, g'=id_{W'}$ we get  $T^{-1}(\mathcal{F}^{-1}_{U,V}\sten'id_{W'} ){\mathcal{F}'}_{T(U\sten V),W'}={\mathcal{F}'}_{{U'\sten V',W'}}(T^{-1}(\mathcal{F}^{-1}_{U,V})\sten'id_{T^{-1}(W')})={\mathcal{F}'}_{U'\sten V',W'}({\mathcal{F}'}_{U',V'}\sten'id_{T^{-1}(W')}$. We call these two relations, auxiliary relations. Now from (2.6) we deduce $\Phi'_{U',V',W'}(id_{U'}\sten'\mathcal{F}^{-1}_{V,W})\mathcal{F}^{-1}_{U,V\sten W}=(\mathcal{F}^{-1}_{U,V}\sten id_{W'})\mathcal{F}^{-1}_{U\sten V,W}T(\Phi_{U,V,W})$ and if we apply $T^{-1}$ to this relation and use the auxiliary relations then we get the counterpart of (2.6).
 
We need to check how the above notions of tensor functor and gauge equivalence extend to $\bar\CC$. Thus, given $(T,\CF)$, we construct a family of isomorphisms $S_\al=S_\al(U_1,\ld ,U_m):(T(U_1)\sten'\cd\sten' T(U_m);\al)\flsh T((U_1\sten\cd\sten U_m;\al)),~\all\al\in\Lan_m,~U_i\in  \mathcal{C}$, by induction on $m$; for $m=1$ we set $S_\al(U)=id_{T(U)}$ and for $m=2$ we set $S_\al(U_1,U_2)=\mathcal{F}_{U_1,U_2}$. Now let  $S_\al$ have been defined already for $\al\in \Lan_k,~k<m$ and let $\gam\in \Lan_m$. Then by definition of the set $\Lan_m$ there exist unique integers $k,l<m$ and unique $\al\in\Lan_k,~\beta\in\Lan_l$ such that $k+l=m$ and $\gam=(\al,\beta)$. Then we set   \eqn{}{S_{(\al,\beta)}(U_1,\ld ,U_m):=\mathcal{F}_{(U_1\sten\cd\sten U_k;\al),(U_{k+1}\sten\cd\sten U_m;\beta)}(S_\al(U_1,\ld ,U_k)\sten'S_\beta(U_{k+1},\ld ,U_m))}
Now let $I'^\al_\beta$ denote the corresponding higher degree associators for $ \mathcal{C}'$, then we claim that 
\eqn{}{S_\beta^{-1}T(I^\al_\beta(U_1,\ld,U_m))S_\al=I'^\al_\beta(U_1',\ld,U_m')}, where $U_i'=T(U_i)$. To prove this, let us denote the left hand side of the above relation by $J^\al_\beta(U_1',\ld,U_m')$. We check that the relations (2.1)-(2.3) hold for this family and therefore by uniqueness of higher degree associators the claim will be proven. We have
$$J^\beta_\gam J^\al_\beta=S_\gam^{-1}T(I^\beta_\gam)S_\beta S_\beta^{-1}T(I^\al_\beta)S_\al=S_\gam^{-1}T(I^\beta_\gam I^\al_\beta)S_\al=S_\gam^{-1}T(I^\al_\gam )S_\al=J^\al_\gam$$ and
\begin{eqnarray}J^{\al_1}_{\al_2}\sten' J^{\beta_1}_{\beta_2}&=&S_{\al_2}^{-1}T(I^{\al_1}_{\al_2})S_{\al_1}\sten' S_{\beta_2}^{-1}T(I^{\beta_1}_{\beta_2})S_{\beta_1}\nonumber\\&=&(S_{\al_2}^{-1}\sten'S_{\beta_2}^{-1})(T(I^{\al_1}_{\al_2})\sten'T(I^{\beta_1}_{\beta_2}))(S_{\al_1}\sten' S_{\beta_1})\nonumber\\&=&S_{(\al_2,\beta_2)}^{-1}\mathcal{F}_{X_2,Y_2}(T(I^{\al_1}_{\al_2})\sten'T(I^{\beta_1}_{\beta_2}))\mathcal{F}_{X_1,Y_1}^{-1}S_{(\al_1,\beta_1)}\nonumber\\&=&S_{(\al_2,\beta_2)}^{-1}T(I^{\al_1}_{\al_2}\sten I^{\beta_1}_{\beta_2})S_{(\al_1,\beta_1)}\nonumber\\&=&S_{(\al_2,\beta_2)}^{-1}T(I^{(\al_1,\beta_1)}_{(\al_2,\beta_2)})S_{(\al_1,\beta_1)}\nonumber\\&=&J^{(\al_1,\beta_1)}_{(\al_2,\beta_2)}\nonumber\end{eqnarray}
where for simplicity we wrote $X_i=(U_1\sten\cd\sten U_k;\al_i),~Y_i=(U_{k+1}\sten\cd\sten U_m;\beta_i)$. And we have
\begin{eqnarray}J^{(id,\sten')}_{(\sten',id)}({U'}_1,{U'}_2,{U'}_3)&=&S_{(\sten',id)}^{-1}T(I^{(id,\sten)}_{(\sten,id)}(U_1,U_2,U_3))S_{(id,\sten')}\nonumber\\&=&(\mathcal{F}_{U_1,U_2}^{-1}\sten id_{T(U_3)})\mathcal{F}_{U_1\sten U_2,U_3}^{-1}T(\Phi_{U_1,U_2,U_3})\mathcal{F}_{U_1,U_2\sten U_3}(id_{T(U_1)}\sten\mathcal{F}_{U_2,U_3})\nonumber\\&=&{\Phi'}_{{U'}_1,{U'}_2,{U'}_3}\nonumber\end{eqnarray}
where we used $S_{(id,\sten')}(U_1,U_2,U_3)=\mathcal{F}_{U_1,U_2\sten U_3}(id_{T(U_1)}\sten\mathcal{F}_{U_2,U_3})$ which comes from (2.9), and similarly for $S_{(\sten',id)}$.\\

Now let $\bar{\mathcal{C}}$ and  $\bar{\mathcal{C}'}$ be the canonical extensions of $\mathcal{C}$ and $\mathcal{C}'$ to strict braided monoidal Ab-categories respectively. We define a functor $\bar{T}:\bar{\mathcal{C}}\flsh\bar{\mathcal{C}'}$ by $ \bar{T}(U_1\sten\cd\sten U_m):=T(U_1)\sten'\cd\sten' T(U_m)$ and for morphism $f:U_1\sten\cd\sten U_m\flsh V_1\sten\cd\sten V_n$ in $ \bar{\mathcal{C}}$ we define $\bar{T}(f)$ as follows; let $f^\al_\beta:(U_1\sten\cd\sten U_m;\al)\flsh(V_1\sten\cd\sten V_n;\beta)$ be the representative of $f$ for any $\al\in\Lan_m,\beta\in\Lan_n$. We define $\bar{T}(f)$ be the equivalence class of the morphism $S_\beta^{-1}T(f^\al_\beta)S_\al:(T(U_1)\sten'\cd\sten' T(U_m);\al)\flsh(T(V_1)\sten'\cd\sten' T(V_n);\beta)$. We must show that this definition does not depend to the choice of the representative. Thus let $f^{\al_i}_{\beta_i},~ i=1,2$ be two representatives of $f$ then we have 
\begin{eqnarray}{I'}^{\beta_1}_{\beta_2}S_{\beta_1}^{-1}T(f^{\al_1}_{\beta_1})S_{\al_1}{I'}^{\al_2}_{\al_1}&=&S_{\beta_2}^{-1}T(I^{\beta_1}_{\beta_2})S_{\beta_1}S_{\beta_1}^{-1}T(f^{\al_1}_{\beta_1})S_{\al_1}S_{\al_1}^{-1}T(I^{\al_2}_{\al_1})S_{\al_2}\nonumber\\&=&S_{\beta_2}^{-1}T(I^{\beta_1}_{\beta_2}f^{\al_1}_{\beta_1}I^{\al_2}_{\al_1})S_{\al_2}\nonumber\\&=&S_{\beta_2}^{-1}T(f^{\al_2}_{\beta_2})S_{\al_2}\nonumber\end{eqnarray}
Thus $\bar{T}(f)$ is well-defined.

Next, for composable morphisms $f$ and $g$ in $\bar{\mathcal{C}}$ we have
\begin{eqnarray}\bar{T}(gf)&=&[S_\gam^{-1}g^\beta_\gam f^\al_\beta S_\al]\nonumber\\&=&[S_\gam^{-1}g^\beta_\gam S_\beta S_\beta^{-1}f^\al_\beta S_\al]\nonumber\\&=&\bar{T}(g)\bar{T}(f)\nonumber
\end{eqnarray}
and for addable morphisms $f$ and $g$ in $\bar{\mathcal{C}}$ we have  $\bar{T}(f+g)=[S_\beta^{-1}(f^\al_\beta+g^\al_\beta)S_\al]=\bar{T}(f)+\bar{T}(g)$. And clearly for all objects $X=U_1\sten\cd\sten U_m,Y=V_1\sten\cd\sten V_n$ in   $\bar{\mathcal{C}}$ we have $\bar{T}(X\sten Y)=\bar{T}(X)\sten'\bar{T}(Y)$ and if $X'=U'_1\sten\cd\sten U'_{m'},Y'=V'_1\sten\cd\sten V'_{n'}$ be another objects in $\bar{\mathcal{C}}$ and $f:X\flsh Y,~g:X'\flsh Y'$ be morphisms in  $\bar{\mathcal{C}}$ we have
\begin{eqnarray}\bar{T}(f\sten g)&=&[S_{(\beta,\beta')}^{-1}T(f^\al_\beta\sten g^{\al'}_{\beta'})S_{(\al,\al')}]\nonumber\\&=&[(S_\beta^{-1}\sten'S_{\beta'}^{-1})\mathcal{F}^{-1}_{Y_\beta,Y'_{\beta'}}T(f^\al_\beta\sten g^{\al'}_{\beta'})\mathcal{F}^{-1}_{X_\al,X'_{\al'}}(S_\al\sten' S_{\al'})]\nonumber\\&=&[(S_\beta^{-1}\sten'S_{\beta'}^{-1})(T(f^\al_\beta)\sten'T( g^{\al'}_{\beta'}))(S_\al\sten' S_{\al'})]\nonumber\\&=&[S_\beta^{-1}T(f^\al_\beta)S_\al\sten'S_{\beta'}^{-1}T( g^{\al'}_{\beta'})S_{\al'}]\nonumber\\&=&\bar{T}(f)\sten'\bar{T}(g)\nonumber\end{eqnarray}
where $X_\al=(U_1\sten\cd\sten U_m;\al),X'_{\al'}=(U'_1\sten\cd\sten U'_m;\al'),Y_\beta=(V_1\sten\cd\sten V_n;\beta),Y'_{\beta'}=(V'_1\sten\cd\sten V'_n;\beta')$.

Now, with with above notation for $X,Y,X_\al$ and $Y_\beta$, we show that 
\begin{eqnarray}\bar{T}(\bar{\Psi}_{X,Y})=\bar{\Psi'}_{\bar{T}(X),\bar{T}(Y)}\end{eqnarray}
Let us denote $\bar{\Psi}_{X,Y}$ by $f$ for brevity. Then by definition for $\al\in\Lan_m,\beta\in\Lan_n$ we have $f^{(\al,\beta)}_{(\beta,\al)}=\Psi_{X_\al,Y_\beta}$. Thus
\begin{eqnarray}\bar{T}(\bar{\Psi}_{X,Y})&=&[S_{(\beta,\al)}^{-1}T(f^{(\al,\beta)}_{(\beta,\al)})S_{(\al,\beta)}]\nonumber\\&=&[(S_\beta^{-1}\sten'S_\al^{-1})\mathcal{F}^{-1}_{Y_\beta,X_\al}T(\Psi_{X_\al,Y_\beta})\mathcal{F}^{-1}_{X_\al,Y_\beta}(S_\al\sten' S_\beta)]\nonumber\\&=&[(S_\beta^{-1}\sten'S_\al^{-1}){\Psi'}_{T(X_\al),T(Y_\beta)}(S_\al\sten' S_\beta)]\nonumber\\&=&[{\Psi'}_{{\bar{T}(X)}_\al,{\bar{T}(Y)}_\beta}]\nonumber\\&=&\Bar{\Psi'}_{\bar{T}(X),\bar{T}(Y)}\nonumber
\end{eqnarray}
where by ${\bar{T}(X)}_\al$ we mean $(T(U_1)\sten'\cd\sten'T(U_m);\al)$ and by ${\bar{T}(Y)}_\beta$ we mean $(T(V_1)\sten'\cd\sten'T(V_n);\beta)$ and in fourth equality we used naturality of braid in $\mathcal{C}'$.

Summarizing the above argument, we have
\begin{prop}\label{barT}
Any tensor functor $(T,\mathcal{F})$ from $(\mathcal{C},\sten,\Phi,\Psi)$ to $(\mathcal{C}',\sten',\Phi',\Psi')$ induces  a canonical additive functor $\bar{T}:\bar{\mathcal{C}}\flsh\bar{\mathcal{C}'}$ obeying

(i)~$\bar{T}(X\sten Y)=\bar{T}(X)\sten'\bar{T}(Y)$, for all objects X,Y in   $\bar{\mathcal{C}}$

(ii)~$\bar{T}(f\sten g)=\bar{T}(f)\sten'\bar{T}(g)$, for all morphisms f,g in   $\bar{\mathcal{C}}$

(iii)~$\bar{T}(\bar{\Psi}_{X,Y})=\Bar{\Psi'}_{\bar{T}(X),\bar{T}(Y)}$, for all objects X,Y in   $\bar{\mathcal{C}}$. 
\end{prop}

\section{Braided Hochschild and Cyclic Cohomology}

Before starting this section let us agree that if in a diagram all the bands are labeled by the object $A$, we label the bands by integers $0,1,2,3,\ldots $. For examples we use the diagram (b) instead of diagram (a) in Figure~2.
\begin{figure}
\[ \epsfbox{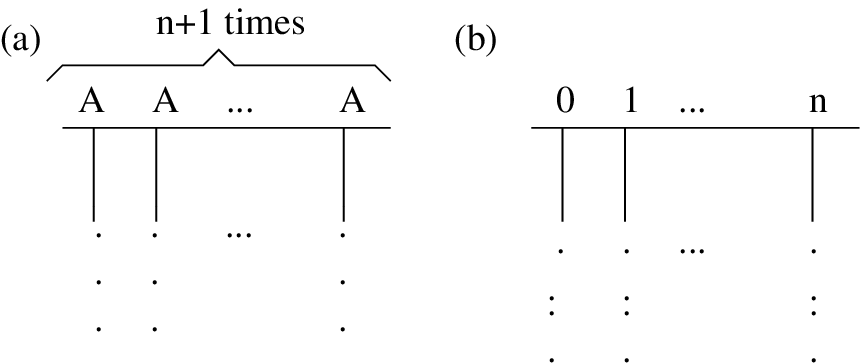} \]
\caption{}
\end{figure} 
We let $\mathcal{C}$ be a braided unital monoidal Ab-category. 
\begin{defi}\label{ribalg} A ribbon algebra in $\mathcal{C}$ is an algebra $(A,m,\eta)$ (see the previous section) equipped with an automorphism $\sig:A\flsh A$ such that
\begin{eqnarray}m(\sig\sten\sig)\Psi^2=\sig m,~~~\sig\eta=\eta\nonumber
\end{eqnarray}
where $\Psi=\Psi_{A,A}$ is the braid $A\sten A\flsh A\sten A$ .
\end{defi}
 If the category  $\mathcal{C}$ is the category of vector spaces over a field with trivial braiding i.e. the flip, then these relations mean that $\sig$ is just an algebra automorphism preserving the unit. In fact this Definition~\ref{ribalg} is a combination of the axioms of an algebra homomorphism and the axiom for the relation between the braiding and the ribbon structure in a ribbon category. It will be rather essential for us in the construction of braided Hochschild and cyclic cohomology. We call  $\sig$ a {\em ribbon automorphism \/} for the algebra $A$. As was mentioned in the previous section, we will work via the strict extension $\bar{\mathcal{C}}$ of $\mathcal{C}$. 

We set $C_n=A^{\sten(n+1)},~~n\ge 0$ as an object of $\bar{\mathcal{C}}$  of degree $n+1$ and define morphisms 
\begin{eqnarray}d_i=d^{(n)}_i:C_n\flsh C_{n-1},~s_i=s_i^{(n)}:C_n\flsh C_{n+1}~,0\le i\le n,~\lan=\lan_n: C_n\flsh C_n 
\end{eqnarray}
in $\bar{\mathcal{C}}$ by diagrams (a) ,(b), (c) and (g) in Figure~3, where all bands in all diagrams in Figure~3 are labeled by object $A$ in $\mathcal{C}$.  We use a diagrammatic notation for morphisms as explained in the Preliminaries.

\begin{figure}
\[\epsfbox{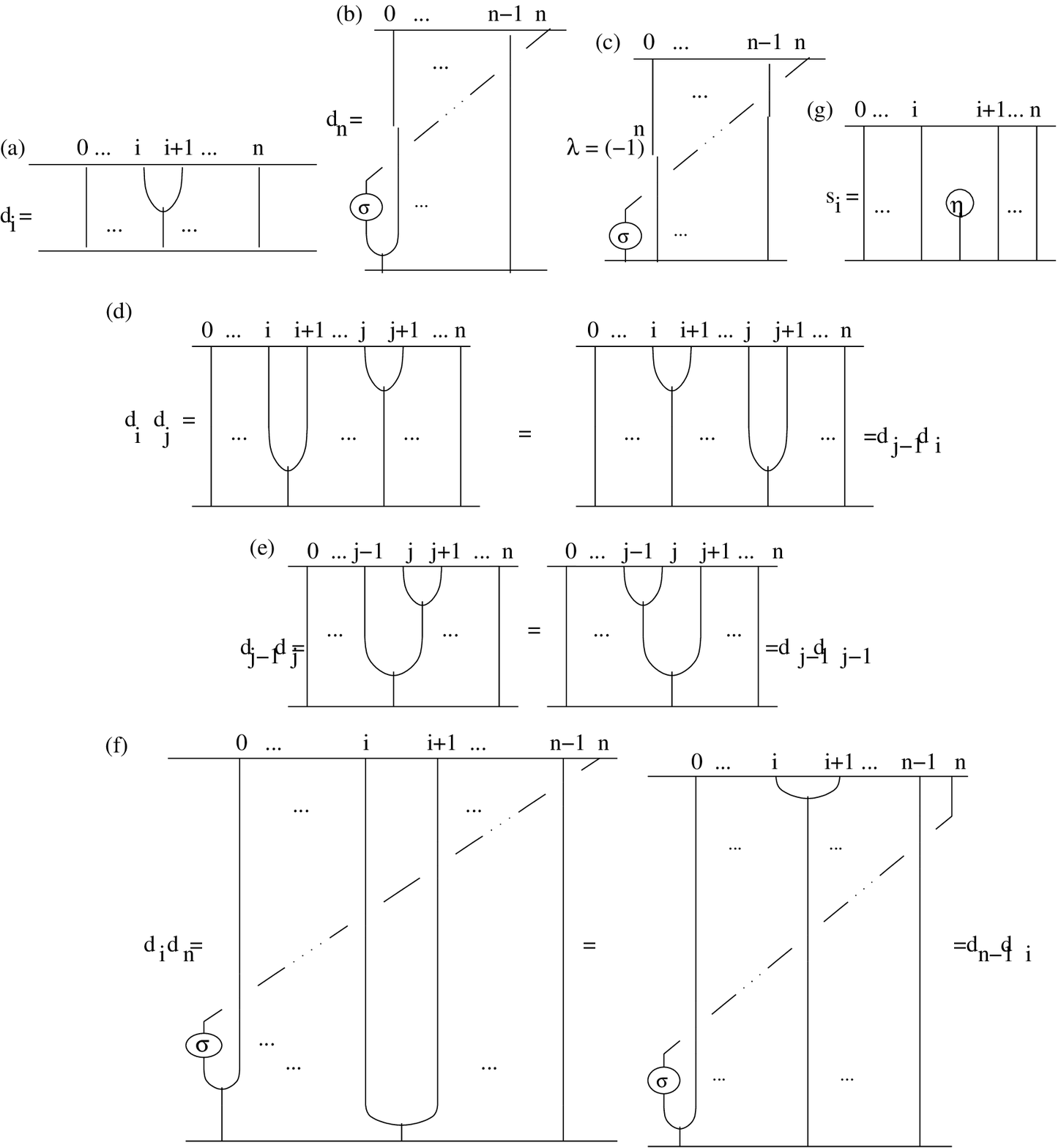}\]
\caption{}
\end{figure}

\begin{thm}\label{simplicial} On $C_n$ we have
$$(i)~d_id_j=d_{j-1}d_i,~ 0\le i< j\le n~~~(ii)~s_is_j=s_{j+1}s_i,~ 0\le i\le j\le n$$ 
$$(iii)~ d_is_j= \left\{ \begin{array}{r@{\quad,\quad}l}s_{j-1}d_i & i<j\\id~~~~ & i=j~ or~i=j+1\\s_jd_{i-1} & i>j+1\end{array}\right.$$
$$(iv)~d_i\lan=-\lan d_{i-1},~1\le i\le n~~~d_0\lan=(-1)^nd_n$$
$$(v)~s_i\lan=-\lan s_{i-1},1\le i\le n~~~s_0\lan=(-1)^n\lan^2s_n$$
$$(vi)~d_i\lan^{n+1}=\lan^nd_i,~~~(vii)~s_i\lan^{n+1}=\lan^{n+2}s_i,~~~0\le i\le n$$      
\end{thm}
\proof  (i)  The proof by means of graphical calculus is in Figure~3 parts (d) (for $j<n,i<j-1$), (e) (for $j<n,i=j-1$), (f) (for $j=n, i<n-1$ ) and Figure~4 part (a) (for $j=n, i=n-1$ ). Note that in the third equality of Figure~4(a) we used the identity in Figure~4(b) and in the fifth equality we used  Definition~\ref{ribalg}.

The proof of parts (ii)-(v) are very straightforward and we leave them to the reader. (vi) At first note that $\lan$ is invertible with inverse given in Figure~5. Now from the recursive relations in first part of (iv) it is easy to compute all the $d_i$ in terms of $d_n$ and $\lan$  and then using the last part of (iv) we compute $d_i$ in terms of $d_0$ and $\lan$ as  $d_i=(-1)^i\lan^{-n+i}d_0\lan^{n-i+1}~~~,0\le i\le n$. In particular $d_0=\lan ^{-n}d_0\lan ^{n+1}$. Thus $d_i=(-1)^i\lan^id_0\lan^{-i}~~~,\all i$. Hence we have
\begin{eqnarray}
d_i\lan^{n+1}&=&(-1)^i\lan ^id_0\lan^{n+1-i}\nonumber\\
&=&(-1)^i\lan^i\lan^nd_0\lan^{-i}\nonumber\\
&=&(-1)^i\lan^n\lan^id_0\lan^{-i}\nonumber\\
&=&\lan^nd_i.\nonumber
\end{eqnarray}
Part (vii) is similar to (iv).
\eproof

\begin{figure}
\[ \epsfbox{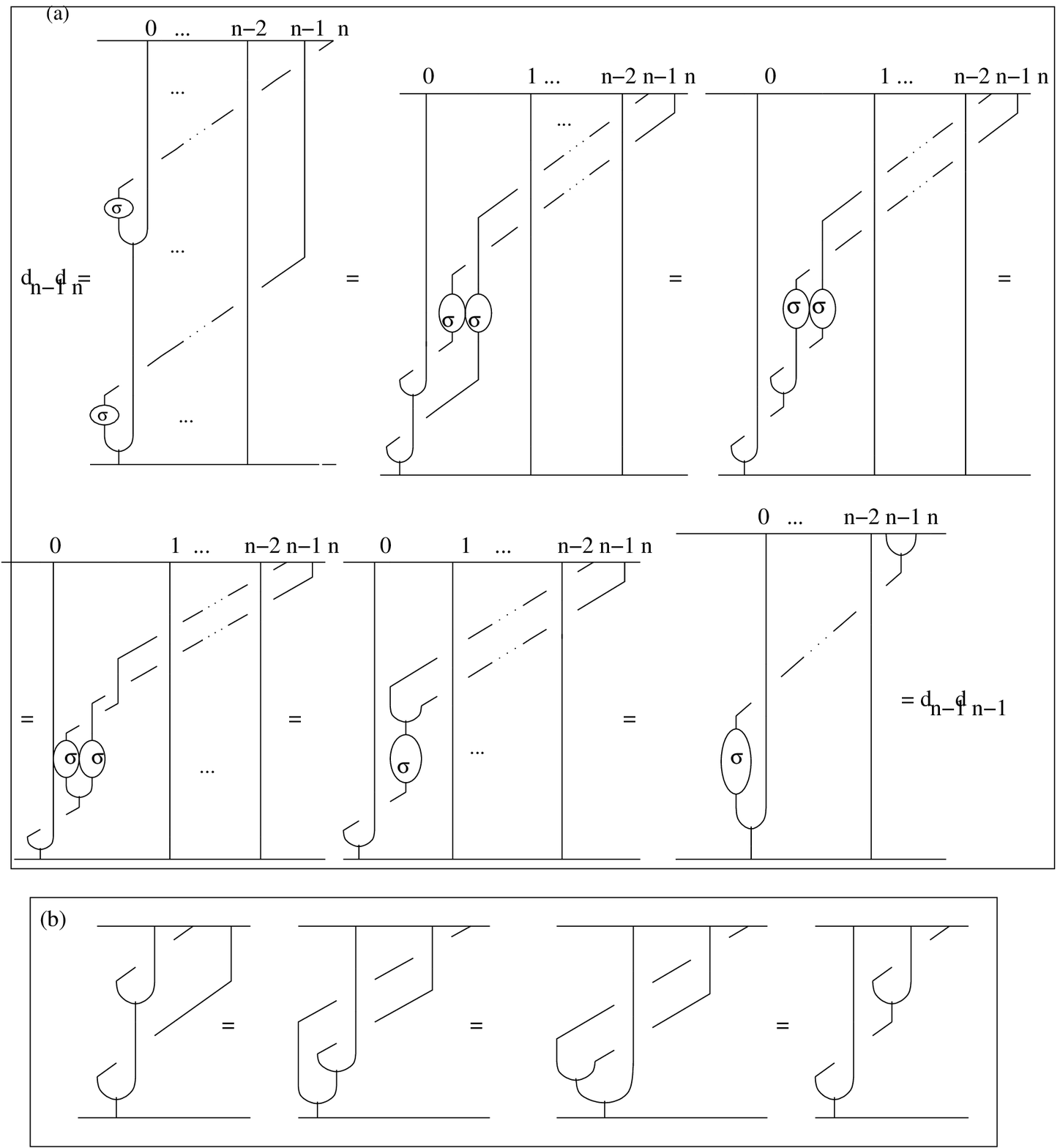} \]
\caption{}
\end{figure}

\begin{figure}
\[ \epsfbox{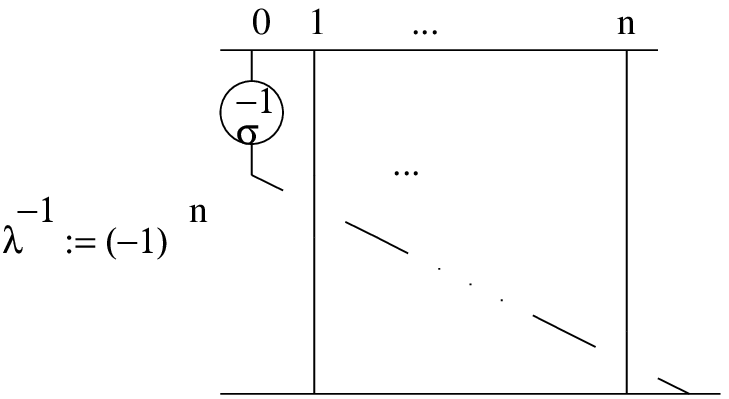} \]
\caption{}
\end{figure}
Following the strategy in \cite{KMT}, we now define 
\eqn{}{C^n(\mathcal{C};A,\sig)=\{\sfi\in Hom_{\bar{\mathcal{C}}}(A^{\sten(n+1)},\mbox{\boldmath $1$})~|~\sfi\lan^{n+1}=\sfi\}}  
By the above Theorem~\ref{simplicial} the morphisms $d_i,\lan$ and $s_i$ induce morphisms
\eqn{}{d_i:C^{n-1}\flsh C^n,~\lan:C^n\flsh C^n,~s_i:C^{n+1}\flsh C^n,~0\le i\le n}
respectively, where we use same symbols. Hence for example $d_i(\sfi):=\sfi d_i,~\sfi\in C^{n-1}$. Then we obtain a cocyclic module $\{C^n\}_{n\ge 0}$ with the above linear maps as face, cyclicity and degeneracy maps respectively. Namely we have    
$$On ~C^{n-1};~~~~~~~~~~~d_jd_i=d_id_{j-1},~~~0\le i<j\le n$$
$$On ~C^{n+1};~~~~~~~~~~~s_js_i=s_is_{j+1},~~~0\le i\le j\le n$$
$$On ~C^n;~~~~~~s_j d_i= \left\{ \begin{array}{r@{\quad,\quad}l}d_is_{j-1} & i<j\\id~~~~ & i=j~ or~i=j+1\\d_{i-1}s_j & i>j+1\end{array}\right.$$

$$On ~C^{n-1};~~~~~\lan d_i=-d_{i-1}\lan,~~~1\le i\le n~~~\lan d_0=(-1)^nd_n$$
$$On~C^{n+1};~~~\lan s_i=- s_{i-1}\lan,~~~1\le i\le n~~~\lan s_0=(-1)^ns_n\lan^2$$ 
$$On ~C^n;~~~~~~~~~\lan^{n+1}=id.$$
Therefore the general theory of Hochschild and cyclic cohomology\cite{Con} gives us a cochain complex $(C^*,d)$,~~$d:=\sum_{i=0}^n(-1)^id_i$  and we call the cohomology of this complex the \emph{ braided Hochschild cohomolgy} of the ribbon  algebra $(A,\sig)$ in the category $\mathcal{C}$ and denote it by $HH^*(\mathcal{C};A,\sig)$. And also we have a subcomplex of the above complex defined as usual by
$$C^n_\lan(\mathcal{C};A,\sig)=\{\sfi\in Hom_{\bar{\mathcal{C}}}(A^{\sten(n+1)},\bf{1})~|~\lan(\sfi)=\sfi\}$$ and we call its cohomology the \emph{ braided cyclic cohomology } of the ribbon algebra $(A,\sig)$ in the category $\mathcal{C}$ and denote it by $HC^*(\mathcal{C};A,\sig)$.

Now we suppose that the ring $K=Hom(\bf{1},\bf{1})$ is a field containing $\mathbb{Q}$, the rational numbers. Then again the general theory of Hochschild and cyclic cohomology gives the SIB-long exact sequence
$$\xymatrix@1{\cdots\ar[r]& HC^{n-1}\ar[r]^{\mathcal{S}}& HC^{n+1}\ar[r]^{\mathcal{I}}& HH^{n+1}\ar[r]^{\mathcal{B}}& HC^n\ar[r]&\cdots }$$
where $\mathcal{I}$ is induced from the inclusion map $C^*_\lan\hookrightarrow C^*$  and $\mathcal{B}$ is implemented by the Connes' boundary map $\mathcal{B}:C^{n+1}\flsh C^n$ defined by $\mathcal{B}=(-1)^nN(s_{-1}+s_n)$, where $N=\sum_{i=0}^n\lan^i$ and $s_{-1}=(-1)^ns_0\lan^{-1}$ is the extra degeneracy map. Finally, $\mathcal{S}$ is the periodicity map which  (see for example \cite{GVF}, formula (10.15)), is given explicitly by \eqn{}{\mathcal{S}(\sfi)=\frac{-1}{n(n+1)}\sum_{1\le i\le j\le n}(-1)^{i+j}d_{i-1}(d_{j-1}(\sfi)).}

Let us now extend the notion of an ordinary  'differential calculus' (DC)  over an ordinary  algebra  (i.e. a DC over an algebra in the category of vector spaces) to an algebra $A$ inside a monoidal Ab-category. Until now we have seen only the universal calculus treated in such generality\cite{Ma:diag}, using diagrammatic methods. In fact the axioms are the same as usual, namely part of a differential graded algebra\cite{Con} including $A$ in degree zero; but this time all objects and morphisms must be inside the category. Thus, a \emph{ differential calculus  of degree $1\le n\le\infty$  over algebra $(A,m,\eta)$, in $\mathcal{C}$} is a sequence of objects $\Om=\{\Om_i\}_{i=0}^n$ in $\mathcal{C}$ together with morphisms $\{m_{i,j}:\Om_i\sten \Om_j\flsh\Om_{i+j}\}_{i,j=0}^n$, called multiplication and morphisms $d=\{d_i:\Om_i\flsh\Om_{i+1}\}_{i=0}^n$, called exterior differentials, such that $\Om_0=A$ , $m_{0,0}=m$, all with the well-known axioms for an ordinary DC when viewed in the category  $\bar{\mathcal{C}}$. For instance the diagram of the Leibniz rule is in Figure 6 part (c) where by labels $i$ and $j$ we mean $\Om_i$ and $\Om_j$. We do need to say some words about the following axiom of an ordinary  DC, where one usually requires that every $k$-form is a sum of  $k$-forms of the form $a^0da^1 \cdots da^k$ for some $a^i\in A$. We translate this axiom for a DC in a category by requiring that for each $k$ the morphisms 
\begin{eqnarray}m(id\sten d): A\sten\Om_{k-1}\flsh\Om_k 
\end{eqnarray}    
be epimorphisms. One can easily conclude by using induction on $k$, that the morphisms
$$m(id_A\sten d^{\sten k}):A^{\sten(k+1)}\flsh\Om_k$$
are epimorphisms of $\bar{\mathcal{C}}$ (the converse is also true but we do not need it).
We recall that a morphism $f:U\flsh V$ in a category is called an {\em epimorphism} if for each pair of  morphisms $g,~h:V\flsh W$ the equality $gf=hf$ implies $g=h$. We suppose in what follows  that the tensor product of two epimorphisms in $\mathcal{C}$  is also an epimorphism.
 
\begin{defi}\label{ribcalc} A ribbon DC  over a ribbon algebra $(A,\sig)$ is a DC,  $\Om$, together with a sequence  of  automorphisms   $\{\sig_i:\Om_i\flsh\Om_i\}_{i=0}^n$  such that $\sig_0=\sig$ and  
\begin{eqnarray}
m_{i,j}(\sig_i\sten\sig_j)\Psi_{j,i}\Psi_{i,j}=\sig_{i+j} m_{i,j},~~~d_i\sig_i=\sig_i d_i
\end{eqnarray}
for all $i,j$, where by $\Psi_{i,j}$ we mean $\Psi_{\Om_i,\Om_j}.$ 
 \end{defi}
This has been presented in Figure 6 part (b). It means that $\sigma$ extends to a ribbon structure on $\Omega$ in a manner compatible with $d$.  
\begin{defi}\label{ribtrace}
A {\em ribbon graded trace} (r.g.t.) on a  degree $n$ ribbon DC is a morphism $\int :\Om_n\flsh \mbox{\boldmath $1$}$  such that if $i+j=n$ then
\begin{eqnarray}\label{rgt} 
\int m_{i,j}=(-1)^{ij}\int m_{j,i}(\sig_j\sten id_{\Om_i})\Psi_{i,j}.
\end{eqnarray} 
 It is called {\em closed} if $\int d=0$
\end{defi}
 This has been presented in Figure~6 part (a). If $\int$ satisfies  (\ref{rgt}) at least for $i=0 ,j=n$  on a (not necessarily ribbon) DC on a ribbon algebra in $\mathcal{C}$ then we call it a \emph{weak ribbon graded trace} (w.r.g.t.).
 
\begin{figure}
\[\epsfbox{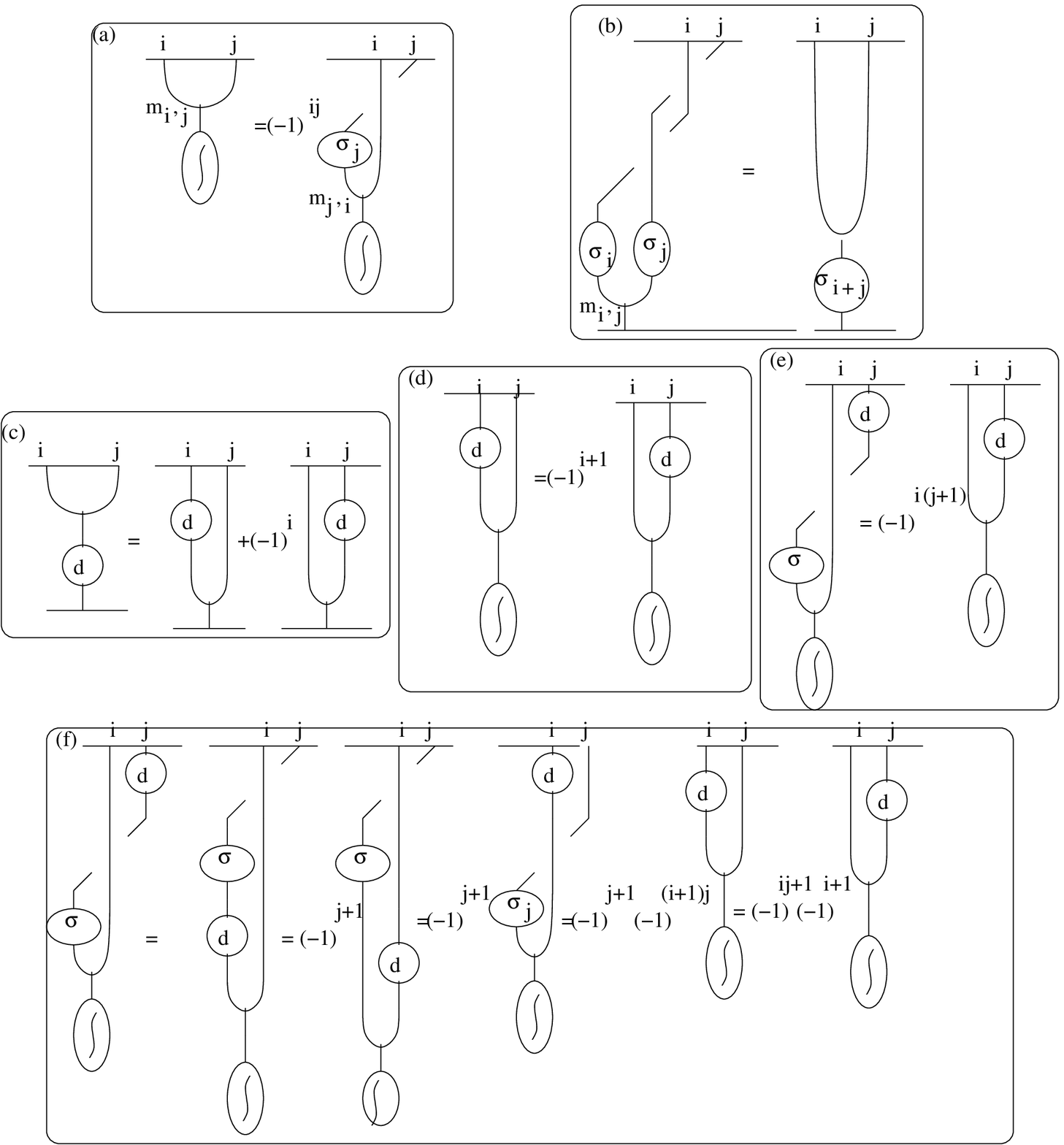}\]
\caption{}
\end{figure}

Following the strategy in \cite{KMT} we have
\begin{prop}\label{wbrgt} 
On a ribbon DC, any w.r.g.t. is also a r.g.t. 
\end{prop}
\proof  We use induction on $j$. Let (\ref{rgt}) be true for $j$; we prove it for $j+1$. We have the identity in Figure~6 part (e) as proven in part (f). In part (f) we used part (d) which is an immediate consequence of the Leibniz rule (represented by the diagrams in part (c)) and closedness of $\int$. In the fourth equality of part (f) we used the induction hypothesis. Now we have the  identity in Figure~7  part (a) as proven in part (b), where $i$ means $\Om_i$, $j$ means $\Om_j$ and $0$ means $\Om_0=A$. Here in the third equality we used the equality in Figure~6 part (e) and in the sixth equality we used (\ref{rgt}). Now since the morphism $m(id\sten d):A\sten \Om_j\flsh\Om_{j+1}$ is an epimorphism, then by our assumption, $id\sten( m(id\sten d)):\Om_i\sten A\sten \Om_j\flsh\Om_i\sten\Om_{j+1}$, is also an epimorphism. Thus from the equality in Figure~7 part (a)  we conclude that the induction hypothesis holds for $j+1$.   
\eproof 
\begin{figure}
\[\epsfbox{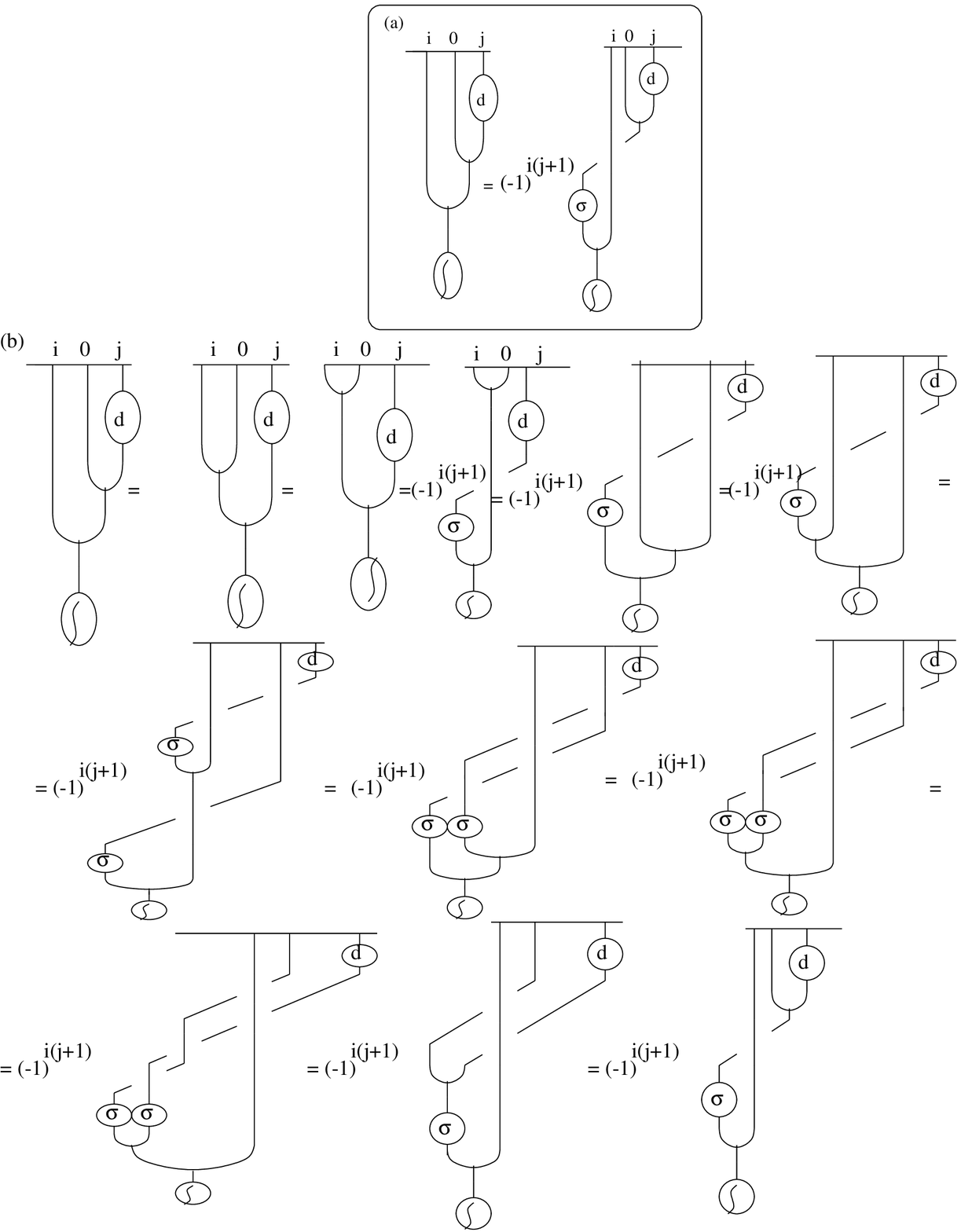}\]
\caption{}
\end{figure}

\begin{thm}\label{character}
Let $\Om$ be a (not necessarily ribbon) DC of degree $0\le n\le\infty$ on the ribbon algebra $(A, \sig)$ in $\mathcal{C}$, and $\int$ a closed w.r.g.t. on $\Om$. Define the morphism in  $\bar{\mathcal{C}}$
\begin{eqnarray}
\sfi:A^{\sten(n+1)}\flsh \mbox{\boldmath $1$},~~~\sfi=\int m(id\sten d^{\sten n})
\end{eqnarray}
(the diagram of $\sfi$ is in Figure~9 (c)). Then $\sfi$ is a braided cyclic cocycle, i.e. \(\sfi\in Z^n_\lan(\mathcal{C} ;A,\sig)\). Here $m$ is the morphism $A\sten\Om_1^{\sten n}\flsh \Om_n$
induced by the multiplication morphisms $m_i$ (using associativity)\end{thm}  
\proof  Let us denote the morphism $md^{\sten k}$ by $f_k$. Then clearly we have $f_k=m(d\sten f_{k-1})=m(f_{k-1}\sten d)$ and $df_k=0$ and by definition we have $\sfi=\int m(id\sten f_n)$. Then the proof that $\lan(\sfi)=\sfi$ is in Figure~8, where the left hand side diagram is by definition $\lan(\sfi)=\sfi$ and by $[k]$ we mean $A^{\sten k}$ .

\begin{figure}
\[\epsfbox{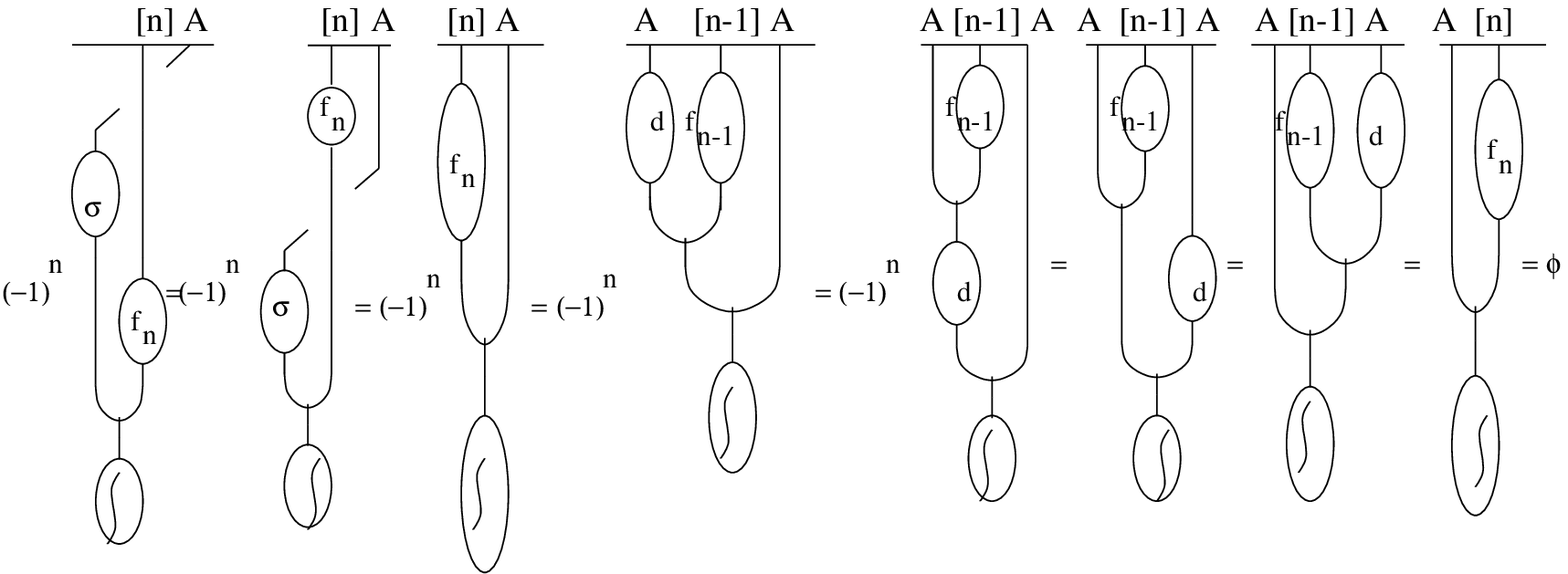}\]
\caption{}
\end{figure}   

\begin{figure}
\[\epsfbox{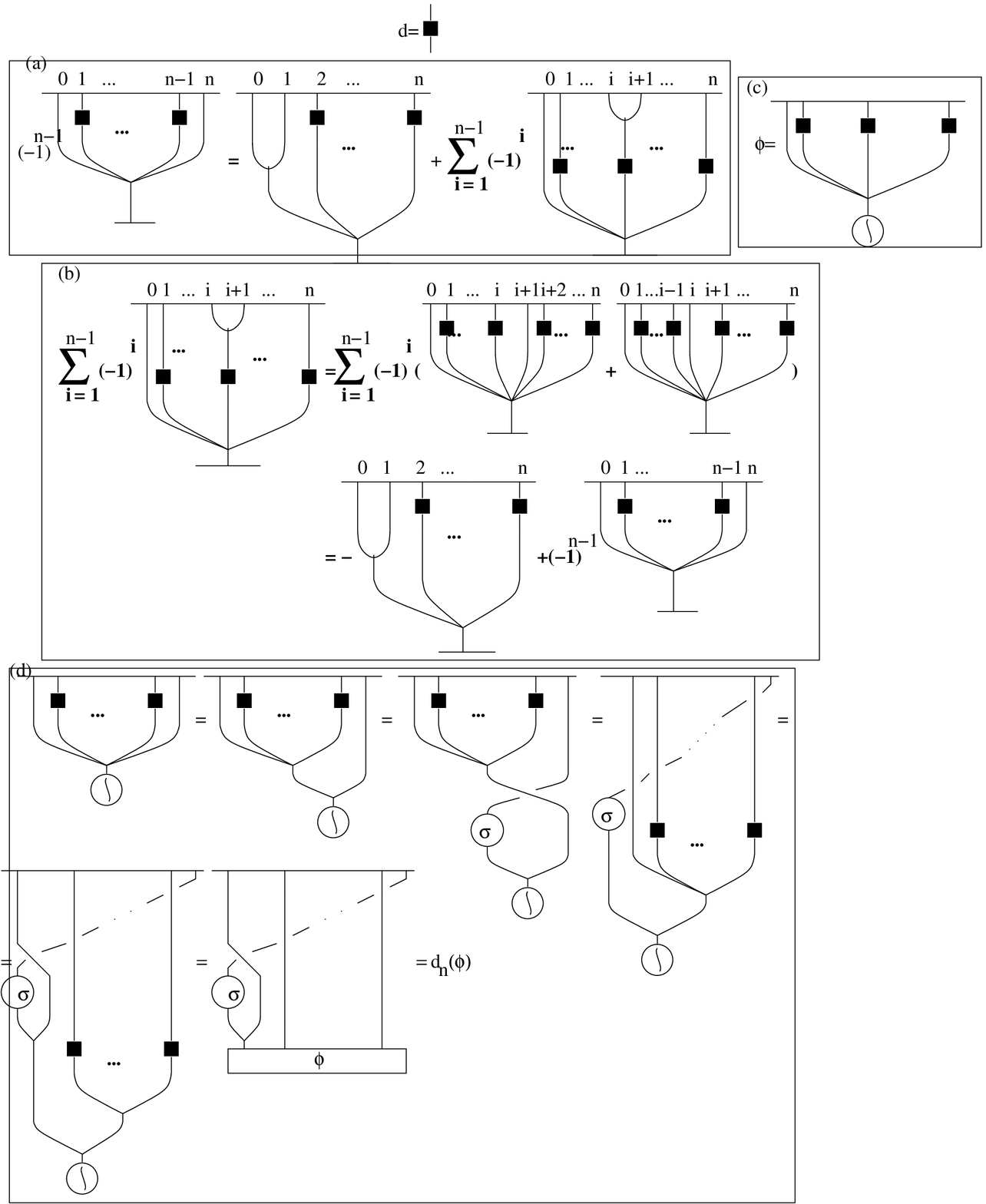}\]
\caption{}
\end{figure}

To prove that $\sfi$ is cocycle we first note the identity in Figure~9 part (a) as proven in part (b). Now applying $\int$ to both sides of this identity and using the definition of $\sfi$, we just need to prove that the left hand side is $d_n(\sfi)$, which is in part (d).
\eproof

Finally, we study the behavior of braided cyclic cohomology with respect to gauge transformation between braided monoidal Ab-categories, using the definitions from the  previous section. 
\begin{thm}\label{gauge}
Let  $(T,\mathcal{F})$ be a tensor functor between braided monoidal Ab-categories $(\mathcal{C},\sten,\Phi,\Psi)$ and $(\mathcal{C}',\sten',\Phi',\Psi')$. And let $(A,m,\eta,\sig)$ be a ribbon algebra in $\mathcal{C}$. Then $A':=T(A)$ is an algebra in  $\mathcal{C}'$ with product $m':=T(m)$,  unit $\eta':=T(\eta)$ and ribbon structure $\sig':=T(\sig)$. Moreover, there exists a morphism of cocyclic modules
$$\bar{T}:C^n(\mathcal{C};A,\sig)\flsh C^n(\mathcal{C}';A',\sig')$$
If in addition $\mathcal{C}$ and $\mathcal{C}'$ are gauge equivalent then the above morphism is a cocyclic module isomorphism and therefore induces an isomorphism between Hochschild and cyclic cohomologies of $A$ and $A'$ in $\mathcal{C}$ and $\mathcal{C}'$ respectively.
\end{thm}
\proof At first suppose $\mathcal{C},\mathcal{C}'$ are strict categories and $\mathcal{F}$ is trivial i.e. $T(U\sten V)=T(U)\sten'T(V)$ and $\mathcal{F}_{U,V}=id$ for all objects $U,V$ in $\mathcal{C}$.  Then from (\ref{tf}) and (\ref{tf-psi}) we get $T(f\sten g)=T(f)\sten'T(g)$ for all morphisms $f,g$ in $\mathcal{C}$, and $T(\Psi_{U,V})={\Psi'}_{T(U),T(V)}$. In this case it is obvious that $(A',m',\eta',\sig')$ is a ribbon algebra in $\mathcal{C}'$. Now if $d_i,s_i,\lan_i$ and $d_i',s_i',\lan_i'$ are the face, degeneracy and cyclicity maps for $A$ and $A'$ respectively then since these maps are constructed by composition, addition or tensor product of the  product of $A$, identity, braiding or ribbon morphisms and since $T$ preserves composition, addition and tensor product of morphisms, we deduce easily that $d_i'=T(d_i),s_i'=T(s_i)$ and $\lan_i'=T(\lan_i)$. Thus the theorem is proved in this case. Now since braided   Hochschild and cyclic cohomology is defined inside the extended strict category $\bar{\mathcal{C}}$, the general non-strict case follows using Proposition~\ref{barT}.
\eproof     
   
\section{Braided cyclic cohomology  of quasialgebras over coquasitringular coquasibialgebras }

In this section we use 'gauge transformation' to construct nontrivial quasialgebras following the
methods in \cite{AlbMa}, and see how the differential calculi and the braided cyclic cohomology behave in this case. We start with the general theory before specializing to the group algebra case  of particular interest. Concrete examples then follow in Section~5.

\subsection{General construction by Drinfeld cotwists} 

We recall\cite{Ma:tan,Ma:book}, cf. \cite{Drinfeld}; 
\begin{defi}
A coquasitriangular coquasibialgebra is  a coalgebra $(H,\Del,\epsilon)$ equipped with a linear map $.:H\sten H\flsh H$, called product, an associated unit element, a convolution invertible ``unital 3-cocycle'', in the sense of a linear map   $ \phi:H^{\sten3}\flsh\mathbb{C}$ satisfing 
$$\phi(b_{(1)},c_{(1)},d_{(1)})\phi(a_{(1)},b_{(2)}.c_{(2)},d_{(2)})\phi(a_{(2)},b_{(3)},c_{(3)})=\phi(a_{(1)},b_{(1)},c_{(1)}.d_{(1)})\phi(a_{(2)}.b_{(2)},c_{(2)},d_{(2)}),$$
$$\phi(a,1,b)=\epsilon(a)\epsilon(b)$$
such that 
$$a_{(1)}.(b_{(1)}.c_{(1)})\phi(a_{(2)},b_{(2)},c_{(2)})=\phi(a_{(1)},b_{(1)},c_{(1)})(a_{(2)}.b_{(2)}).c_{(2)}$$
$\all a,b,c,d\in H$, and finally a convolution invertible linear map $ \mathcal{R}:H\sten H\flsh\mathbb{C} $, satisfing 
$$\mathcal{R}(a.b,c)=\phi( c_{(1)},a_{(1)},b_{(1)})\mathcal{R}(a_{(2)},c_{(2)})\phi^{-1}(a_{(3)},c_{(3)},b_{(2)})\mathcal{R}(b_{(3)},c_{(4)})\phi(a_{(4)},b_{(4)},c_{(5)}),$$
$$\mathcal{R}(a,b.c)=\phi^{-1}(b_{(1)},c_{(1)},a_{(1)})\mathcal{R}(a_{(2)},c_{(2)})\phi(b_{(2)},a_{(3)},c_{(3)})\mathcal{R}(a_{(4)},b_{(3)})\phi^{-1}(a_{(5)},b_{(4)},c_{(4)}),$$
$$b_{(1)}.a_{(1)}\mathcal{R}(a_{(2)},c_{(2)})=\mathcal{R}(a_{(1)},b_{(1)})a_{(2)}.b_{(2)}$$ 
\end{defi}

Also we recall that the category of (left) $H$-comodules (abbreviated as $H$-\emph{Com} ) is  a braided monoidal Ab-category with 
$$\Phi_{U,V,W}:U\sten (V\sten W)\flsh(U\sten V)\sten W,~~~u\sten v\sten w\mapsto \phi(u_{(1)},v_{(1)},w_{(1)})u_{(2)}\sten v_{(2)}\sten w_{(2)}$$
$$\Psi_{U,V}:U\sten V\flsh V\sten U,~~~u\sten v\mapsto\mathcal{R}(v_{(1)},u_{(1)})v_{(2)}\sten u_{(2)}$$
$\all ~U ,V,W\in\emph{H-Com}$.  An algebra $A$ in this category is called a \emph{left quantum quasispace} (or  left H-comodule quasialgebra\cite{AlbMa}). This means that  we have a left $H$-comodule structure, $A\flsh H\sten A, ~a\mapsto a_{(1)}\sten a_{(2)}$, a linear map $.:A\sten A\flsh A$, called product, which is unital, associative in the category and equivariant under the coaction of $H$. Specifically, 
$$a.(b.c)= \phi(a_{(1)},b_{(1)},c_{(1)})(a_{(2)}.b_{(2)}).c_{(2)},~~~(a.b)_{(1)}\sten (a.b)_{(2)}=a_{(1)}.b_{(1)}\sten a_{(2)}.b_{(2)},~a,b\in A.$$ Finally, a ribbon structure for $A$ in the sense of Section~3 means a $H$-comodule isomorphism $\sig:A\flsh A$ such that
\eqn{sig}{\sig(a.b)=\mathcal{R}(b_{(1)}, a_{(1)})\mathcal{R}(a_{(2)}, b_{(2)})\sig(a_{(3)}).\sig(b_{(3)})}

Now let us describe explicitly a braided cyclic cocycle in this category. First of all by the argument in Section~2,  for $H$- comodules $U_i,~1\le i\le n$, a morphism   $U_1\sten\cdots\sten U_n\flsh \mathbb{C}$ in the category of $H$-Com, is an equivalence class  of $H$-comodule intertwiners $U_1\sten_\al\cdots \sten_\al U_n\flsh \mathbb{C}$ for all $\al\in\Lan_n$. Here a $H$-comodule intertwiner $U_1\sten_\al\cdots \sten_\al U_n\flsh \mathbb{C}$ means a linear map $f:U_1\mbox{\boldmath{$\sten$}} \cdots\mbox{\boldmath{$\sten$}} U_n\flsh \mathbb{C}$ satisfying
\eqn{mor}{f(u^1\sten\cd\sten u^n)1=u^1_{(1)}._\al\cdots ._\al u^n_{(1)}f(u^1_{(2)}\sten\cd\sten u^n_{(2)}),~~~\all u_i\in U_i,}
where in order to avoid confusion, we have used the boldmath notation  $U_1\mbox{\boldmath{$\sten$}} \cdots\mbox{\boldmath{$\sten$}} U_n$ for the usuall vector space tensor product, $u_1\sten\cd\sten u_n$ meaning an element of the vector space $U_1\mbox{\boldmath{$\sten$}} \cdots\mbox{\boldmath{$\sten$}} U_n$ and the notation $u^1._\al\cdots ._\al u^n$ as defined in Section~2. Similarly, a morphism $U_1\sten\cdots\sten U_m\flsh V_1\sten\cdots\sten V_n$ is an equivalence class of $H$-comodule intertwiners $U_1\sten_\al\cdots \sten_\al U_m\flsh V_1\sten_\beta\cdots \sten_\beta V_n,~~~\al\in\Lan_m,\beta\in\Lan_n$.

Now let us describe the morphisms $d_i:C_n=A^{\sten(n+1)}\flsh C_{n-1}=A^{\sten n},~0\le i\le n-1$ in our particular category. These are represented by morphisms
\eqn{b_i}{A\sten_\al\cdots\sten_\al(A\sten A)\sten_\al\cdots\sten_\al A\flsh A\sten_\al\cdots\sten_\al A\nonumber}
\eqn{}{(a^0\sten\cd\sten a^n)\mapsto (a^0\sten\cd\sten a^i.a^{i+1}\sten\cd\sten a^n),~~~a^i\in A ,\al\in\Lan_n \nonumber}
while $d_n$ by definition is $(m_A\sten id_{n-1})(\sig_A\sten id_n)\Psi_{n,1}$, where $ \Psi_{n,1}:=\Psi_{A^{\sten n},A}$ and $id_k:=id_{A^{\sten k }}.$
 Let us write for simplicity $X$ for  $A^{\sten^{n-1}_\al} ,\al\in\Lan_{n-1}$. Then the representative of $d_n$ from $(A\sten X)\sten A$ to $A\sten X$ is $(m_A\sten id_X)\phi_{A,A,X}(\sig \sten id_{A\sten X})\Psi_{A\sten X,A}$ which has value on $a^0\sten\cdots\sten a^n$ equal to
\begin{eqnarray}&&\kern-30pt \mathcal{R}(a^n_{(1)},a^0_{(1)}.(a^1_{(1)}._\al\cdots._\al a^{n-1}_{(1)})) (m_A\sten id_X)\phi_{A,A,X}(\sig \sten id_{A\sten X})(a^n_{(2)}\sten(a^0_{(2)}\sten(a^1_{(2)}\sten_\al\cdots\sten_\al a^{n-1}_{(2)})))\nonumber\\
&=&\mathcal{R}(a^n_{(1)},a^0_{(1)}.(a^1_{(1)}._\al\cdots._\al a^{n-1}_{(1)})) (m_A\sten id_X)\phi_{A,A,X}(\sig(a^n_{(2)})\sten(a^0_{(2)}\sten(a^1_{(2)}\sten_\al\cdots\sten_\al a^{n-1}_{(2)}))\nonumber\\
&=&\mathcal{R}(a^n_{(1)},a^0_{(1)}.(a^1_{(1)}._\al\cdots._\al a^{n-1}_{(1)}))\phi(a^n_{(2)},a^0_{(2)},a^1_{(2)}._\al\cdots._\al a^{n-1}_{(2)})(m_A\sten id_X)((\sig(a^n_{(3)})\sten a^0_{(3)})\sten\nonumber\\
&&(a^1_{(3)}\sten_\al\cdots\sten_\al a^{n-1}_{(3)}))\nonumber\\
&=&\mathcal{R}(a^n_{(1)},a^0_{(1)}.(a^1_{(1)}._\al\cdots._\al a^{n-1}_{(1)}))\phi(a^n_{(2)},a^0_{(2)},a^1_{(2)}._\al\cdots._\al a^{n-1}_{(2)})\sig(a^n_{(3)}). a^0_{(3)}\sten a^1_{(3)}\sten\cdots\sten a^{n-1}_{(3)}\nonumber
\end{eqnarray}
Therefore the representative of $d_n$ from $(A\sten(A^{\sten^{n-1}_\al}))\sten A$ to $ A\sten(A^{\sten^{n-1}_\al})$ is
\begin{eqnarray}
d_n(a^0\sten\cdots\sten a^n)&=&\mathcal{R}(a^n_{(1)},a^0_{(1)}.(a^1_{(1)}._\al\cdots._\al a^{n-1}_{(1)}))\phi(a^n_{(2)},a^0_{(2)},a^1_{(2)}._\al\cdots._\al a^{n-1}_{(2)})\times\nonumber\\
& &\sig(a^n_{(3)}). a^0_{(3)}\sten a^1_{(3)}\sten\cdots\sten a^{n-1}_{(3)},
\end{eqnarray}
where $\al\in\Lan_{n-1}$.

We can similarly represent $\lan$ with the representative $(A\sten(A^{\sten^{n-2}_\al}))\sten A\flsh A\sten(A\sten(A^{\sten^{n-1}_\al}))$ by 
\eqn{}{\lan(a^0\sten\cdots\sten a^n)=\mathcal{R}(a^{n-1}_{(1)},a^0_{(1)}.(a^1_{(1)}._\al\cdots._\al a^{n-2}
_{(1)}))\sig(a^{n-1}_{(2)})\sten a^0_{(2)}\sten\cdots\sten a^{n-2}_{(2)},}
where $\al\in\Lan_{n-2}$.

Next a DC,~ $\Om(A)=\bigoplus_{i=0}^n\Om_i(A)$,  on $A$ in this category is just a left $H$-covariant DC as in \cite{Wor}, except that the product of forms are associative up to the associator, namely
$$\om.({\om}'.{\om''})=\phi(\om_{(1)},{\om'}_{(1)},{\om''}_{(1)})(\om_{(2)}.{\om'}_{(2)}).{\om''}_{(2)},~~\all \om,{\om'},{\om''}\in\Om.$$
A w.r.g.t.  on $\Om(A)$ is a left $H$-invariant functional $\Om_n\flsh \mathbb{C}$ satisfying
\eqn{}{\int\om.a=\mathcal{R}( a_{(1)},\om_{(1)})\int\sig(a_{(2)})\om_{(2)}.}

Now we recall briefly the so called Drinfeld  cotwist and associated ``gauge transformation''. For this discussion we will denote the product of two elements $a$ and $b$ in $H$ or $A$, just by $ab$ and we keep the notation $a\ast b$ or  $a.b$ respectively for the new products defined as follows. Thus for $(H,\Phi,\mathcal{R})$ as above and any convolution invertible linear map $F:H\sten H\flsh \mathbb{C}$ obeying $F(a,1)=F(1,a)=\epsilon(a),~\all a\in H$ (a 2-cochain), we define a coquasitringular coquasibialgebra $H^F$ as follows; $H^F$  as a coalgebra is, by definition, $(H, \Del,\epsilon)$ itself and the product for $H^F$ is defined by \cite{Drinfeld}, 
\eqn{}{a\ast b:= F(a_{(1)},b_{(1)})a_{(2)}b_{(2)}F^{-1}(a_{(3)},b_{(3)}),}
 the unital 3-cocycle for $H^F$ is defined by
\eqn{phi-F}{\phi_F(a,b,c):=F(b_{(1)},c_{(1)})F(a_{(1)},b_{(2)}c_{(2)})\phi(a_{(2)},b_{(3)},c_{(3)})F^{-1}(a_{(3)}b_{(4)},c_{(4)})F^{-1}(a_{(4)},b_{(5)})}
and the coquasitringular structure for $H^F$ is defined by
\eqn{psi-F}{\mathcal{R}_F(a,b):=F(b_{(1)},a_{(1)})\mathcal{R}(a_{(2)},b_{(2)})F^{-1}(a_{(3)},b_{(3)})}
$\all a,b\in H$. This is the dual version of the  Drinfeld twist, (see \cite{Ma:book}). Then, as above, we have a braided monoidal Ab-category, namely $H^F$-Com, which is nothing other  than $H$-Com, as an Ab-category but with new monoidal and braided structures. We denote the tensor product, the associator and the braiding of this category by $\sten_F,\Phi_F,\psi-F$ respectively.

If $A$ is a quantum quasispace over $H$ then there is a gauge transformed copy of it in $H^F$-Com, namely we cotwist the product of $A$ as \cite{AlbMa} 
$$a.b:=F(a_{(1)},b_{(1)})a_{(2)}b_{(2)},~~~\all a,b\in A$$
then $A_F$ is by definition $A$ as a left $H$-comodule equipped with the  above product. It is easy to see that $A_F$ is a (left) quantum  quasispace  over $H^F$, see \cite{AlbMa}.            

Now let us show that if $\sig$ is a ribbon structure for $A$ in $H$-Com, then  $\sig$ is still a ribbon structure for $A_F$   in the category of left $H^F$-Com. Indeed,
\begin{eqnarray}  
&&\kern -40pt \mathcal{R}_F(b_{(1)}, a_{(1)}))\mathcal{R}_F(a_{(2)}, b_{(2)})\sig(a_{(3)}).\sig(b_{(3)})\nonumber\\
&=&F(a_{(1)}, b_{(1)})\mathcal{R}(b_{(2)}, a_{(2)})F^{-1}(b_{(3)}, a_{(3)})F(b_{(4)}, a_{(4)}) \mathcal{R}(a_{(5)}, b_{(5)})F^{-1}(a_{(6)}, b_{(6)})\nonumber\\
&&\quad\quad F(a_{(7)}, b_{(7)})\sig(a_{(8)})\sig(b_{(8)})\nonumber\\
&=&F(a_{(1)}, b_{(1)})\mathcal{R}(b_{(2)}, a_{(2)})\mathcal{R}(a_{(3)}, b_{(3)})\sig(a_{(4)})\sig(b_{(4)})\nonumber\\
&=&F(a_{(1)}, b_{(1)})\sig(a_{(2)}b_{(2)})\nonumber\\
&=&\sig(a. b)\nonumber
\end{eqnarray}
as required.

Next it is known that the cotwisted comodule algebra $\Om(A_F):=\Om(A)_F$ is a DC over the algebra $A_F$ in the category of $H^F$-Com, (see \cite{Ma:cli}). We recall that differential forms in  $\Om(A_F)$ are as the same as differential forms in $\Om(A)$ with same coaction of $H$ on them and the differential operator, $d$, remain unchanged but the product of forms has been cotwisted to 
$$\om.\om':=F(\om_{(1)},{\om'}_{(1)})\om_{(1)}{\om'}_{(2)}$$
Therefore for left invariant forms this product does not change, i.e. we have
$$ \om.\om'=\om\om',~~~a.\om=a\om,~~~\om.a=\om a,~~~\all\om\in\Om^{inv},~a\in A$$
Where  $\Om^{inv}$ denote the space of left invariant forms.

Let us show now that if $\int$ is a w.r.g.t. on $\Om(A)$ then it is also a w.r.g.t. on $\Om(A_F)$. Indeed, 
\begin{eqnarray}
\mathcal{R}_F(a_{(1)},\om_{(1)})\int\sig(a_{(2)}).\om_{(2)}&=& F(\om_{(1)}, a_{(1)})\mathcal{R}( a_{(2)},\om_{(2)})F^{-1}(a_{(3)}, \om_{(3)})
F(a_{(4)}, \om_{(4)})\int\sig(a_{(5)})\om_{(5)}\nonumber\\
&=&\int F(\om_{(1)}, a_{(1)})\om_{(2)} a_{(2)}\nonumber\\
&=&\int\om.a\nonumber
\end{eqnarray}
as required. 

\emph{Remark}. Observe that even if the associator and coquasitringular structure are trivial before gauge transformation, they are typically no longer trivial after gauge transformation. Even if $F$ is a Hopf algebra 2-cocycle\cite{Ma:book}, the associator after gauge transformation becomes trivial but the coquasitiangular structure is not typically trivial but is cotriangular. And note that the cochain $F$ needed for gauge transformation of classical semisimple Lie groups  to their quantum counterpart are not Hopf algebra 2-cocycle but are 2-cochains, thus these observations imply that we do need to consider the class of braided   cyclic cocycles during such transformations.  

As for any monoidal category $\mathcal{C}$ the functor $\mathcal{C}^n\flsh\mathcal{C},~~(U_1,\ld,U_n)\mapsto (\cd((U_1\sten U_2)\sten U_3)\sten\cd)\sten U_n$ is called the \emph{left}-\emph{to}-\emph{right arranger}.   We can use this to give explicit left-to-right representatives of all formulae. Finally,  the categories $(H\emph{-Com},\sten,\Phi,\Psi)$ and $(H^F\emph{-Com},\sten_F,\Phi_F,\psi-F)$ are gauge equivalent, which is the categorical meaning of the Drinfeld cotwist, see \cite{Ma:tan}. We  choose $T=id$ and $\mathcal{F}_{U,V}(u\sten v):=F(u_{(1)},v_{(1)})u_{(2)}\sten v_{(2)}$. Here $\mathcal{F}_{U,V}$ is a morphism in $H^F$-Com 
since
\begin{eqnarray*}&&(u\sten v)_{(1)}\sten \mathcal{F}_{U,V}((u\sten v)_{(2)}=u_{(1)}\ast v_{(1)}\sten F(u_{(2)},v_{(2)})u_{(3)}\sten v_{(3)}\\
&&\quad =F(u_{(1)},v_{(1)})u_{(2)}v_{(2)}F^{-1}(u_{(3)},v_{(3)})\sten  F(u_{(4)},v_{(4)})u_{(5)}\sten v_{(5)}\\
&&\quad =F(u_{(1)},v_{(1)})u_{(2)}v_{(2)}\sten u_{(3)}\sten v_{(3)}={(\mathcal{F}_{U,V}(u\sten v))}_{(1)}\sten{(\mathcal{F}_{U,V}(u\sten v))}_{(2)}\end{eqnarray*}
Similarly, $\mathcal{F}^{-1}_{U,V}(u\sten v)=F^{-1}(u_{(1)},v_{(1)})u_{(2)}\sten v_{(2)}$ is a morphism in $H^F$-Com. Moreover, (\ref{tf}) is a consequence of the definition of a morphism in the category of $H$-Com and (\ref{tf-phi}), (\ref{tf-psi}) are consequences of (\ref{phi-F}) and (\ref{psi-F}). Therefore by Theorem~\ref{gauge}  we have:\begin{thm}There is an isomorphism of cocyclic modules
$$ \Bar{id}:C^k(H\emph{-Com}; A, \sig)\flsh C^k(H^F\emph{-Com};A_F,\sig),~~~\sfi\mapsto \sfi_F$$ 
For instance  if $\sfi^{lr}$ be the left-to-right representative of $\sfi$ then the left-to-right representative of $\sfi_F$, which we denote it by $\sfi_F^{lr}$ is  \begin{eqnarray}\sfi_F^{lr}(a^0,\ldots,a^k)&=&F(a^0_{(1)},a^1_{(1)})F(a^0_{(2)}a^1_{(2)},a^2_{(1)})F((a^0_{(3)}a^1_{(3)})a^2_{(2)},a^3_{(1)})
\cdots\\
& & F((\cdots((a^0_{(k)}a^1_{(k)})a^2_{(k-1)})\cdots)a^{k-1}_{(2)},a^k_{(1)})
\sfi^{lr}(a^0_{(k+1)},a^1_{(k+1)},a^2_{(k)}\ldots,a^k_{(2)})\nonumber
\end{eqnarray}
\end{thm}   
\proof  This is  Theorem~\ref{gauge} in our present case. In (4.9) we just used the definition of $S_\al$ (see Section~2).
\eproof

\subsection{Braided   cyclic cocycles of group algebras and $G$-graded quasialgebras}  

Let $G$ be a group and $H=\mathbb{C}G$ the group algebra on $G$. This is a Hopf algebra with 
$$\Del g=g\sten g,~~~\eps(g)=1,~~~S(g)=g^{-1},~~~\all g\in G$$   
we recall that a 3-unital cocycle for $H=\mathbb{C}G$ is just the linear extension of a group cocycle $\phi:G^3\flsh\mathcal{C}-\{0\}$ in the sense of a function $\phi:G^3\flsh\mathbb{C}-\{0\}$ obeying
$$\phi(g_1,g_2,g_3)\phi(g_0,g_1g_2,g_3)\phi(g_0,g_1,g_2)=\phi(g_0,g_1,g_2g_3)\phi(g_0g_1,g_2,g_3),~~~\phi(g_0,e,g_1)=1,~\all g_i\in G.$$

Of interest is the case when $G$ is abelian. Since  $H$ is cocommutative and associative, we can consider $\mathbb{C}G$ as a coquasibialgebra with any  $\phi$, including to start with  the  trivial $\phi\equiv 1$. Then in this standard case, since we assume $G$ is abelian, a coquasitriangular structure is just a group bicharacter $\mathcal{R}:G^2\flsh\mathbb{C}-\{0\}$ in the sense that, a function obeying
$$\mathcal{R}(g_0g_1,g_2)=\mathcal{R}(g_0,g_2)\mathcal{R}(g_1,g_2),~~~\mathcal{R}(g_0,g_1g_2)=\mathcal{R}(g_0,g_1)\mathcal{R}(g_0,g_2),~~~\all g_i\in G.$$
Again to start with we can choose the trivial $\mathcal{R}\equiv 1$, so that $\mathbb{C}G$  is considered  as a coquasitriangular coquasibialgebra with trivial associator and trivial coquasitriangular structure.\\

Next a left $\mathbb{C}G$-comodule means precisely a \emph{G-}graded vetor space $V$ with coaction of $\mathbb{C}G$ on it is given by $\al(v)=|v|\sten v$ where $|v|\in G$ denotes the degree of a homogeneous vector $v\in V$\cite{Ma:book}. An algebra $A$ in the category of $\mathbb{C}G$-comodules is just a G-graded algebra (recall that we have chosen trivial associator for $\mathbb{C}G$). Now since for every bialgebra $H$ (that is a coquasibialgebra with trivial associator), $A=H$ is an $H$-comodule algebra with the coproduct of $H$ taken as coaction of $H$ on $A$, we can consider the  algebra $A=\mathbb{C}G$ in the category of $\mathbb{C}G$-comodules. Since we have chosen  the trivial coquasitriangular structure for $H=\mathbb{C}G$, a ribbon structure $\sig:A=\mathbb{C}G\flsh A=\mathbb{C}G$ is just a group homomorphism $\sig :G\flsh G$ extended by linearity to $\mathbb{C}G$\\

Finally, let $\{\chi_i:G\flsh \mathbb{C}-\{0\}\}_{i=1}^n$ be a finite set of group characters. We extend each $\chi_i$ on $\mathbb{C}G$ by linearity and denote the extended map  still with $\chi_i$. Clearly we have $\chi_i(ab)=\chi_i(a)\chi_i(b)$~~~$\all a,b\in A,\all i$. It is well-known that if $\Lan$ is a $n$-dimensional vector space with basis $\{\om_i\}_{i=1}^n$ then there exists a unique left covariant FODC, $\Gamma$, on $A=\C G$ such that $\all g\in G,~\all i$
\eqn{om-eqn}{\om_ig=\chi_i(g)g\om_i,~~~dg=\sum_{i=1}^n(\chi_i(g)-1)g\om_i} 
In fact these calculi are bicovariant and the space of right invariant forms $\Gam^{r.inv}$, coincides with the space of left invariant forms, $\Gam^{l.inv}=\Lan$. Thus the Woronowicz braiding $\Psi:\Gamma\sten_A\Gamma\flsh\Gamma\sten_A\Gamma$  is just the map, $\Psi(a\om\sten_A\om')=a\om'\sten_A\om,~~~\all\om,\om'\in \Gam^{l.inv}$. Hence the relations among basis 1-forms $\om_i$ in the DC, $\Om:=T(\Gamma)/ker(\Psi-id)$,  are $\om_i^2=0,~~~\om_i\om_j=-\om_j\om_i,~~~\all i,j$, where $T(\Gamma)$  is the tensor algebra over $\Gamma$. Therefore we have a top form $\theta:=\om_1\cdots\om_n$ for the space of left invariant $n$-forms, i.e. the space of left invariant $n$-forms is one dimensional. Now let us define $\pi:\Om_n\flsh A$ and $\rho:A\flsh A$ as follows; since $\{\theta\}$ is a free A-basis for $\Om_n$, for each $\om\in\Om_n $ and $a\in A$ there exist unique elements of $A$, $\pi(\om)$ and $\rho(a)$ such that $\om=\pi(\om)\theta$ and $\theta a=\rho(a)\theta$. One can easily verify that both $\pi$ and $\rho$ are morphisms in the category of $H$-comodules, and $\rho$ is an algebra automorphism and by very definition we have $ \pi(a\om)=a\pi(\om),~~~\pi(\om a)=\pi(\om)\rho(a)$. In our case of $A=\mathbb{C}G$ by (\ref{om-eqn}) we have $$\rho(g)=\chi_1(g)\cdots\chi_n(g)g.$$

Next there exists a unique left and right invariant functional on $\mathbb{C}G$ defined 
by\\
\eqn{}{h(g)=0,~~~\all g\neq e,~~~ h(e)=1}
which defines a canonical left invariant functional 
\eqn{int-G}{\int:\Om_n\flsh \mathbb{C},~~~\int\om:=h(\pi(\om)).}
Since $\om a=\pi(\om)\theta a=\pi(\om)\rho(a)\theta=\rho(a)\pi(\om)\theta=\rho(a)\om$, we have $\int\om a=\int\rho(a)\om$. Thus $\int$ is a w.r.g.t. with ribbon morphism  $\sig(g)=\chi_1(g)\cdots\chi_n(g)g$. Let us show that it is closed. Since $\int$ is left covariant we have $(\int  dg_1\cdots dg_n)e=g_1\cdots g_n \int dg_1\cdots dg_n,~~\all g_i\in G$ thus if $g_1\cdots g_n\neq e$ then $\int  dg_1\cdots dg_n=0$, and if  $g_1\cdots g_n= e$ then we have 
\begin{prop}
If for $g_0,\ldots,g_k\in G$ we have $g_0\cdots g_k=e$, then $dg_0\cdots dg_k=0$
\end{prop}
\proof  We have $g^{-1}dg=\sum_{i=1}^n(\chi_i(g)-1)\om_i$ and $(dg^{-1})g=\sum_{i=1}^n(\chi_i(g^{-1})-1)g^{-1}\om_ig=\sum_{i=1}^n(\chi_i(g^{-1})-1)\chi_i(g)\om_i=-\sum_{i=1}^n(\chi_i(g)-1)\om_i$ Thus $dg^{-1}dg=(dg^{-1})gg^{-1}dg=-\sum_{i=1}^n(\chi_i(g)-1)(\chi_j(g)-1)\om_i\om_j=0$. Thus $dg^{-1}dg=0$ and $d\om(g)=d(g^{-1}dg)=dg^{-1}dg=0$. Since every left invariant form $\om\in\Om_{inv}$ is generated by left invariant 1-forms, we deduce that $d\om=0,~~\all\om\in\Om_{inv}$. Now we show by induction on $k$ that for arbitrary $g_1,\ldots, g_k\in G$ we have $(g_1\cdots g_k)^{-1}dg_1\cdots dg_k\in\Om^k_{inv}$; for $k=1$ it is clear, and since for every $g\in G$ by (\ref{om-eqn}) we have $g\Om_{inv}=\Om_{inv}g$ we conclude by using the induction hypothesis that $(g_1\cdots g_k)^{-1}dg_1\cdots dg_k=g_k^{-1}((g_1\cdots g_{k-1})^{-1}dg_1\cdots dg_{k-1})dg_k\in\Om^k_{inv}$. Now since $g_0=(g_1\cdots g_k)^{-1}$ we see that $dg_0\cdots dg_k=d((g_1\cdots g_k)^{-1}dg_1\cdots dg_k)=0$.
\eproof 

Hence the character $\sfi(g_0,\cdots , g_n)= \int  g_0 dg_1\cdots dg_n$ is a trivially braided   cyclic cocycle. Now we calculate $\sfi$ explicitly. First of all, using (\ref{om-eqn}) one has 
\begin{eqnarray} dg_1\cdots dg_n &=& \sum_{i_1,\ldots,i_n =1}^n(\chi_{i_1}(g_1)-1)\cdots(\chi_{i_n}(g_n)-1)\chi_{i_1}(g_2)\cdots\chi_{i_1}(g_n)\chi_{i_2}(g_3)\cdots\chi_{i_2}(g_n)\nonumber\\
&&\quad\quad\quad\quad  \cdots\chi_{i_{n-1}}(g_n)  g_1\cdots g_n\om_{i_1}\cdots\om_{i_n}\nonumber      
\end{eqnarray}
Now for any permutation  $\tau\in S_n$ we have $\om_{\tau(1)}\cdots\om_{\tau(n)}=sgn(\tau)\om_1\cdots\om_n$. Therefore $\sfi(g_0,\cdots, g_n)=0$, if $g_0\cdots g_n\neq e$ and if $g_0\cdots g_n= e$ then
\begin{eqnarray}
\sfi(g_0,\cdots, g_n)&=&\sum_{\tau\in S_n}sgn(\tau)(\chi_{\tau(1)}(g_1)-1)\cdots(\chi_{\tau(n)}(g_n)-1)\chi_{\tau(1)}(g_2)\cdots\chi_{\tau(1)}(g_n)\nonumber\\
& &\quad\quad \chi_{\tau(2))}(g_3)\cdots\chi_{\tau(2)}(g_n)\cdots\chi_{\tau(n-1)}(g_n)\nonumber
\end{eqnarray}
Using the fact that $\chi_i$ are group characters we conclude that
\begin{eqnarray}
\sfi(g_0,\cdots, g_n)
&=&\sum_{\tau\in S_n}sgn(\tau)(\chi_{\tau(1)}(g_1\cdots g_n)-\chi_{\tau(1)}(g_2\cdots g_n))\nonumber\\
& &\quad\quad (\chi_{\tau(2)}(g_2\cdots g_n))-\chi_{\tau(2)}(g_3\cdots g_n))\cdots(\chi_{\tau(n)}(g_n)-1)
\end{eqnarray}
Hence we have  a braided monoidal Ab-category, namely the category of $G$-graded vector spaces with trivial braiding and trivial associator, and we have an algebra $A$ in it i.e. $A=\mathbb{C}G$ itself but with a nontrivial ribbon structure $\sig: G\flsh G$, a group homomorphism given by $\sig(g)=\chi_1(g)\cdots\chi_n(g)g$. Moreover, we have a DC on $A$ and a closed w.r.g.t. and therefore the corresponding   cyclic cocycle.  This is therefore an example of the theory in \cite{KMT}. 

Moreover, since we are working with a strict category we can suppress the subscript $\al$ in the expression $U_1\sten_\al\cdots\sten_\al U_n$ and $a^0._\al\cdots ._\al a^n$.
Therefore the maps $d_i$ and $\lan$ becomes as usual. Let us describe them in this case more explicitly.
First of all a morphism $\sfi:A^{\sten(k+1)}\flsh \mathbb{C}$ in this category is determined by a function $\sfi:G^{k+1}\flsh \mathbb{C}$ such that $\sfi(g_0,\ldots,g_k)e=g_0\cdots g_k\sfi(g_0,\ldots,g_k)$, but since $G$ is a vector space  basis for $\mathbb{C}G$, we conclude  
\eqn{}{Hom_{G-Vec}((\mathbb{C}G)^{\sten(k+1)},\mathbb{C})=\{\sfi:G^{k+1}\flsh \mathbb{C}~|~\sfi(g_0,\ldots,g_k)=0,~~~if~~g_0\cdots g_k\neq e\}}
and
\eqn{}{
(b\sfi)(g_0,\ldots,g_{k+1})=\sum_{i=0}^k(-1)^i\sfi(g_0,\ldots,g_ig_{i+1},\ldots,g_{k+1})+(-1)^{k+1}\chi(g_{k+1})\sfi(g_{k+1}g_0,\ldots,g_k)}
\eqn{}{(\lan\sfi)(g_0,\ldots,g_k)=(-1)^k\chi(g_k)\sfi(g_k,g_0,\ldots,g_{k-1})}
where $\chi=\chi_1\cdots\chi_n$. Note that  since $(\lan^{k+1}\sfi)(g_0,\ldots,g_k)=(-1)^{k(k+1)}\chi(g_0\cdots g_k)\sfi(g_0,\ldots,g_k)$ we conclude that
\begin{eqnarray*}&&\kern -20pt C^k(G-Vec;\mathbb{C}G,\bar\chi)\\
&& =Hom_{G-Vec}((\mathbb{C}G)^{\sten(k+1)},\mathbb{C})=\{\sfi:G^{k+1}\flsh \mathbb{C}~|~\sfi(g_0,\ldots,g_k)=0,~~~if~~g_0\cdots g_k\neq e\},\end{eqnarray*} where $\bar\chi(g):=\chi_1(g)\cdots\chi_n(g)g$.

Now we are ready to cotwist all of the above to obtain a braided cyclic cocycle. Thus, let $F:G^2\flsh \mathbb{C}-\{0\}$ be a function such that $F(g,e)=F(e,g)=1,~~\all g\in G$ and $H=\mathbb{C}G$ as above. Then after cotwisting, $H^F$ has the same product as $H$, because;
$$g_1.g_2=F(g_1,g_2)g_1g_2F(g_1,g_2)^{-1}=g_1g_2,~~~\all g\in G$$
But we have nontrivial 3-cocycle
\eqn{}{\phi_F(g_1,g_2,g_3)=F(g_2,g_3)F(g_1,g_2g_3)F(g_1,g_2)^{-1}F(g_1g_2,g_3)^{-1}}and nontrivial cotriangular structure
\eqn{}{\mathcal{R}_F(g_1,g_2)=F(g_2,g_1)F(g_1,g_2)^{-1}.}

Therefore we have a cotriangular coquasibialgebra which we denote it by $\mathbb{C}^FG$ and an algebra in the category of $\mathbb{C}^FG$-comodules is called a $G$\emph{-graded quasialgebra} \cite{AlbMa} i.e. a $G$-graded algebra $A=\bigoplus_{g\in G}A_g$ such that 
\eqn{}{a(bc)=\phi_F(|a|,|b|,|c|)(ab)c}
on homogeneous elements. Recall that we chose $A=\mathbb{C}G$ as an algebra in the category of $G$-graded vector spaces which cotwists to a G-graded quasialgebra $A_F=\C_FG$ with the product
\eqn{}{g_1.g_2=F(g_1,g_2)g_1g_2}
We also have
\eqn{}{\int\om. g=F(g,\om_{(1)})F(\om_{(1)},g)^{-1}\chi_1(g)\cdots\chi_n(g)\int g.\om_{(2)}}       
while the isomorphism of Theorem~12 becomes
\eqn{sfi-F}{\sfi_F(g_0,\ldots,g_k)=F(g_0,g_1)F(g_0g_1,g_2)\cdots F(g_0\cdots g_{k-1},g_k)\sfi(g_0,\ldots,g_k).}

\section{Examples} 

In this section we collect examples that demonstrate key aspects of the theory above. They are all constructed using cotwisting on coquasitriangular Hopf algebras in Section~4. The first two are
based on abelian groups as in Section~4.2. 

\subsection{Octonions.} Let $G=\mathbb{Z}_2\times\mathbb{Z}_2\times\mathbb{Z}_2$. Thus $\mathbb{C}G$ is an associative algebra with generators $\{1,u,v.w\}$ and relations\\ 
$$u^2=1,~v^2=1,~w^2=1,~uv=vu,~uw=wu,~vw=wv$$
For simplicity we will use the notation $\vec{i}=u^{i_1}v^{i_2}w^{i_3}$ and $\vec{i}+\vec{j}=u^{ i_1+j_1}v^{ i_2+j_2}w^{ i_3+j_3}$ as well. In fact this notation is statement of the  group  $\mathbb{Z}_2^3$ as an additive group.\\
Each $\mathbb{Z}_2$ has a unique nonzero differential calculus (the universal one) and it is natural to equip $G$ with the 'direct product' of three copies of this. To do this, we choose three characters
$$ \chi_1(u^iv^jw^k)=(-1)^i,~~\chi_2(u^iv^jw^k)=(-1)^j,~~\chi_3(u^iv^jw^k)=(-1)^k,~~i,j,k=0,1$$
Then by (\ref{om-eqn}) we have\\ 
$~~~~~~~~~~~~~~~~~~~~~~~~~~~~~~~~~\om_1u=-u\om_1,~~~\om_2u= u\om_2,~~~\om_3u= u\om_3$\\
$~~~~~~~~~~~~~~~~~~~~~~~~~~~~~~~~~~\om_1v=v \om_1,~~~\om_2v=-v\om_2,~~~\om_3v=v \om_3$\\
$~~~~~~~~~~~~~~~~~~~~~~~~~~~~~~~~~\om_1w=w \om_1,~~~\om_2w=w \om_2,~~~\om_3w=-w\om_3$\\
and again by (\ref{om-eqn}) we have\\
$~~~~~~~~~~~~~~~~~~~~~~~~~~~~~~~~~du=-2u\om_1,~dv=-2v\om_2,~dw=-2w\om_3$\\
and\\
$~~~~~~~~~~~~~~~~~~~~~~~~~~~~~~~~~(du)u=-udu,~~(dv)u=udv,~~(dw)u=udw$\\
$~~~~~~~~~~~~~~~~~~~~~~~~~~~~~~~~~(du)v=vdu,~~(dv)v=-vdv,~~(dw)v=vdw$\\
$~~~~~~~~~~~~~~~~~~~~~~~~~~~~~~~~~(du)w=wdu,~~(dv)w=wdv,~~(dw)w=-wdw$\\
and therefore\\
$~~~~~~~~~~~~~~~~~~~~~~~~~~~~~~~~~(du)^2=(dv)^2=(dw)^2=0$\\
$~~~~~~~~~~~~~~~~~~~~~~~~~~~~~~~~~dudv=-dvdu,~~dudw=-dwdu,~~dvdw=-dwdv.$\\
Now using the Leibniz rule and the above relations we have $d(uv)=-2uv(\om_1+\om_2)$,~ $d(uw)=-2uw(\om_1+\om_3)$ and $d(vw)=-2vw(\om_2+\om_3)$ and thus we have $d(uvw)=(du)vw+u(dvw)=-2uvw(\om_1+\om_2+\om_3).$ Therefore, generally, we have
$$ d(u^{i_1}v^{i_2}w^{i_3})=-2u^{i_1}v^{i_2}w^{i_3}(i_1\om_1+i_2\om_2+i_3\om_3),~~~~~i_k=0,1,~k=1,2,3$$ 
or in additive notation
$$d\vec{i}=-2\vec{i}(i_1\om_1+i_2\om_2+i_3\om_3).$$
Using this formula and the above commutation relations we obtain
\begin{eqnarray} \pi(\vec{i}d\vec{j}d\vec{k}d\vec{l})&=&-8(\vec{i}+\vec{j}+\vec{k}+\vec{l})[(-1)^{l_2+ l_3}l_1((-1)^{k_2}j_2 k_3-(-1)^{k_3}j_3 k_2)+\nonumber\\& &(-1)^{l_1+ l_3}l_2((-1)^{k_3}j_3 k_1-(-1)^{k_1}j_1 k_3)+(-1)^{l_1+ l_2}l_3((-1)^{k_1}j_1 k_2-(-1)^{k_2}j_2 k_1)]\nonumber
\end{eqnarray} 
 Thus we calculate the character of the trivially braided ribbon graded trace defined by (\ref{int-G}) as
\begin{eqnarray}\sfi(\vec{i},\vec{j},\vec{k},\vec{l})&=&-8[(-1)^{l_2+ l_3}l_1((-1)^{k_2}j_2 k_3-(-1)^{k_3}j_3 k_2)+(-1)^{l_1+ l_3}l_2((-1)^{k_3}j_3 k_1\nonumber\\
& &-(-1)^{k_1}j_1 k_3)+(-1)^{l_1+ l_2}l_3((-1)^{k_1}j_1 k_2-(-1)^{k_2}j_2 k_1)]
\end{eqnarray}  
if $\vec{i}+\vec{j}+\vec{k}+\vec{l}=0$, and zero otherwise. One can compute this from the general formula (4.13) as well.

Now we study cotwisting of $G=\mathbb{Z}_2^3$ to the octonions. The complex numbers, the quaternions, the octonions and the higher Cayley algebras can be constructed by the Cayley-Dixon process. On the other hand it has been shown in \cite{AlbMa} that these algebras are  $G$-graded quasialgebras of the form $\mathbb{C}_FG$ for $G$ a power of $\mathbb{Z}_2$ and for a suitable $F$.

For octonions we take\cite{AlbMa} $G=\mathbb{Z}_2^3$ and 
\eqn{}{F(\vec{i},\vec{j})=(-1)^{i_1(j_1+j_2+j_3)+i_2(j_2+j_3)+i_3j_3+j_1i_2i_3+i_1j_2i_3+i_1i_2j_3}}
which has coboundary 
\eqn{}{\phi_F(\vec{i},\vec{j},\vec{k})=(-1)^{|\vec{i}~\vec{j}~\vec{k}|}}
where $ |\vec{i}~\vec{j}~\vec{k}|$ is a short notation for determinant of the matrix whose columns are the vectors $\vec{i},\vec{j}$ and $\vec{k}$ respectively. 

We denote the cotwisted product by . i.e, $x.y=F(x,y)xy,~~\all x,y\in \mathbb{C}G$. Then we have the following relations
$$u.u=-1,~v.v=-1,~w.w=-1,~u.v=-v.u,~u.w=-w.u,~v.w=-w.v,~u.(v.w)=-(u.v).w$$          
and the relations between 0-forms and left invariant forms  remain unchanged after twisting the above DC on $G=\mathbb{Z}_2\times\mathbb{Z}_2\times\mathbb{Z}_2$. But we have
$$du.u=-u.du,~~dv.u=-u.dv,~~dw.u=-u.dw\nonumber$$
$$du.v=-v.du,~~dv.v=-v.dv,~~dwv=-v.dw$$
$$du.w=-w.du,~~dv.w=-w.dv,~~dw.w=-w.dw\nonumber$$  
and
$$du.du=dv.dv=dw.dw=0\nonumber$$
$$ du.dv=dv.du,~~du.dw=dw.du,~~dv.dw=dw.dv$$
This is our natural differential calculus or 'exterior algebra' for the octonions as obtained by cotwisting. Like the octonions themselves, it is a nonassociative quasialgebra. We see that now the function algebra generators and their differentials uniformly anticommute, while the latter mutually commute.

From the above, we  have $d(u.v)=du.v+u.dv=-2u.\om_1.v-2u.v.\om_2=-2u.v.(\om_1+\om_2)$ and similarly $d(u.w)=-2u.w.(\om_1+\om_3)$ and $d(v.w)=-2v.w.(\om_2+\om_3)$ and $d((u.v).w)=d(u.v).w+(u.v).dw=-2(u.v.(\om_1+\om_2)).w-2(u.v).w.\om_3=-2(u.v).w.(\om_1+\om_2)-2(u.v).w.\om_3=-2(u.v).w.(\om_1+\om_2+\om_3)$. Therefore, generally, we have
$$d((u^i.v^j).w^k)=-2(u^i.v^j).w^k(i\om_1+j\om_2+k\om_3),~~~~~i,j,k=0,1$$
Using formula (\ref{sfi-F}) we can then calculate 
$$\sfi_F(\vec{i},\vec{j},\vec{k},\vec{l})=F(\vec{i},\vec{j})F(\vec{i}+\vec{j},\vec{k})F(\vec{i}+\vec{j}+\vec{k},\vec{l})\sfi(\vec{i},\vec{j},\vec{k},\vec{l})$$
as the left-to-right representative of the character. 

\subsection{Noncommutative algebraic torus}

Let $G=\mathbb{Z}\times\mathbb{Z}$. Thus $\mathbb{C}G$ is an associative algebra with free commuting generators $u,v$. The standard calculus on algebra $\mathbb{C}G$, i.e. the two dimensional bicovariant DC can be written with basis $\{\om_1,\om_2\}$  for left invariant $1$-forms and relations 
$$ \om_k a=a\om_k,~~~d(u^iv^j)=u^iv^j(i\om_1+j\om_2),~~\om_k^2=0,~~\om_1\om_2=-\om_2\om_1$$
 $\all a\in \mathbb{C}G,\all i,j\in \mathbb{Z}$. This is obtained from characters as in the finite group case via a limiting procedure \cite{Ma:cli}. 

Then we have $\pi(u^{i_1}v^{i_2}d(u^{j_1}v^{j_2})d(u^{k_1}v^{k_2})=(j_1k_2-j_2k_1)u^{i_1+j_1+k_1}v^{i_2+j_2+k_2}$ and hence 
\eqn{}{\sfi(u^{i_1}v^{i_2},u^{j_1}v^{j_2},u^{k_1}v^{k_2})=j_1k_2-j_2k_1}
if $i_1+j_1+k_1=i_2+j_2+k_2=0$, and zero otherwise.

Next, as in \cite{Ma:cli} we chose the cocycle  
$F(u^{i_1}v^{i_2},u^{j_1}v^{j_2})=e^{\imath\theta i_2j_1}$ to gauge transform the category of $\mathbb{Z}\times\mathbb{Z}$-graded vector spaces. Then the algebra $\mathbb{C}_FG$ after cotwisting of the product has the relations $v.u=e^{\imath\theta}u.v$, which  we call \emph{algebraic noncommutative torus}. This observation itself is well-known already for the
full noncommutative torus $C^*$-algebra in a more explicit (non-categorical) context \cite{Rieffel}. 
From our point of view the algebra is associative since above $F$ is a cocycle, but it still gives a nontrivial example of the theory of Section~4. 
 
Now using (4.22) we compute the character of above DC after gauge transformation as
\begin{eqnarray}
\sfi_F(u^{i_1}v^{i_2},u^{j_1}v^{j_2},u^{k_1}v^{k_2})&=&F(u^{i_1}v^{i_2},u^{j_1}v^{j_2})F(u^{i_1}v^{i_2}u^{j_1}v^{j_2},u^{k_1}v^{k_2})\sfi(u^{i_1}v^{i_2},u^{j_1}v^{j_2},u^{k_1}v^{k_2})\nonumber\\
&=&e^{\imath\theta (i_2j_1+(i_2+j_2)k_1)}(j_1k_2-j_2k_1)
\end{eqnarray}
Since, by the very definition of $F$, we have $u^{.i}v^{.j}=u^iv^j, ~\all i,j\in \mathbb{Z}$, we conclude that 
\begin{eqnarray}
\sfi_F(u^{.i_1}v^{.i_2},u^{.j_1}v^{.j_2},u^{.k_1}v^{.k_2})=e^{\imath\theta (i_2j_1+(i_2+j_2)k_1)}(j_1k_2-j_2k_1)
\end{eqnarray}
We remark that we have $v^{.j}u^{.i}=e^{\imath\theta ij}u^iv^j, ~\all i,j\in \mathbb{Z}$. Note that $\sigma$ is trivial in this example. 

\subsection{Quantum quasimanifolds covariant under quantum groups $\C_q(G)$}

Here we take the initial Hopf algebra to be $H=\C(G]$, the algebraic version of a classical Lie group
of complex simple Lie algebra $g$. More precisely, we need to work in a deformation-theoretic setting where $H=\C(G)[[\hbar]]$ is extended over this ring $\C[[\hbar]]$ of formal powerseries in a deformation parameter $\hbar$, rather than working over $\C$ itself. With this proviso, we are able to use the theory in Section~4.1 with $F$ a cochain with values in $\C[[\hbar]]$. By essentially dualising the theory in \cite{Drinfeld}, one knows that there exists an $F$ such that
\[  \C_q(G)\isom (\C(G)[[\hbar]])^F\]
i.e. such that after cotwisting one obtains the (formal power-series version) of the standard quantum group $\C_q(G)$. Here $q=e^{\hbar/2}$. Note that although the required  $F$ is not a cocycle so that the cotwist on the right hand side here is in theory a coquasiHopf algebra, its coboundary $\phi_{KZ}=\partial F$ (which is obtained by solving the Knizhnik-Zamolochikov (KZ) equations) happens to be central in the sense
\[ a b c= \phi_{KZ}(a\o,b\o,c\o) a\t b\t c\t \phi^{-1}_{KZ}(a\thr,b\thr,c\thr)\]
so that the coquasiHopf algebra on the right hand side happens to remain associative as indeed is $\C_q(G)$ on the left hand side (since it is a usual Hopf algebra). This point of view has been expounded recently in \cite{BegMa}. Let us note only that it is not exactly the one of Drinfeld even before dualisation. For that one should start with $\C(G)[[\hbar]]$ as a nontrivial coquasi-Hopf algebra with a nontrivial initial coquasitriangular structure $\CR_0$ built from the Killing form and a certain $\phi_{0}$ as the initial associator, which is the object actually obtained by Drinfeld by solving the KZ-equations. Here $\phi_{0}$ is closely related to $\phi_{KZ}$ above as its `inverse' in the sense that cotwisting it by $F$ gives $\eps$ (the trivial associator). In this way, Drinfeld's theory in the cotwist form would cotwist a certain coquasiHopf algebra $(\C(G)[[\hbar]],\CR_0,\phi_{0})$ into the ordinary Hopf algebra $\C_q(G)$ with its usual quasitriangular structure and trivial associator. By contrast in \cite{BegMa} one starts with $\C(G)[[\hbar]]$ as completely classical with trivial coquasitriangular structure and trivial associator, and obtains after cotwisting $(\C_q(G),\CR_{F},\phi_{KZ})$ as a cotriangular coquasi-Hopf algebra. Like in our examples based on abelian groups, this happens to be an ordinary Hopf algebra in its algebra and coalgebra. Here $\CR_F$ is given by (\ref{psi-F}) and $\phi_{KZ}$ by (\ref{phi-F}). This cotriangular coquasiHopf algebra $(\C_q(G),\CR_{F},\phi_{KZ})$ is the object under which our examples of quantum quasispaces  below are covariant. As an algebra and coalgebra it coincides with the usual quantum group $\C_q(G)$.

Thus, in \cite{BegMa} this  $(\C_q(G),\CR_{F},\phi_{KZ})$ construction was used to obtain a quantum differential calculus on $\C_q(G)$ as
a supercoquasiHopf algebra $\Omega(\C(G))^F$. We extend this setting now to any manifold $M$ on which the classical Lie group $G$ acts. More precisely, we assume that there are algebraic versions
$\C(M)$ for the coordinate algebra and for a coaction $\C(M)\to \C(G)\tens \C(M)$. Thus $\C(M)$ is given as an algebra in our initial category of $\C(G)$-comodules. Moreover, we extend all this data to the formal powerseries setting (we adjoin $\hbar$). 

\begin{thm} Let $M$ be a classical $G$-manifold in the sense above. Then there is a quasialgebra
$\CC_q(M)=(\C(M)[[\hbar])_F$ in the category of $\C_q(G)$-comodules, where $\C_q(G)$ is the standard quantum group associated to $G$ viewed as a cotriangular coquasiHopf algebra. Moreover, $\CC_q(M)$ has a quasiassociative differential calculus $\Omega(\C(M))_F$.
\end{thm}
\proof We apply the theory of Section~4.1 with $H=\C(G)[[\hbar]]$, $F$ the cochain above and $A=\C(M)[[\hbar]]$. \eproof

Similarly, any classical data on $M$ such as a cyclic cocycle twists to a braided cyclic one on the
symmetric but nontrivially monoidal category of $\C_q(G)$-comodules.

Let us note also that when all of our data is obtained from exponentiating infinitesimal
data, we can look at the structure of $\CC_q(M)$ to lowest order. Then one finds
\[\{ \{a,b\},c\}+ \{\{b,c\},a\}+\{\{c,a\},b\}=2\tilde n(a\tens b\tens c)\]
where $n=[r_{+13},r_{+23}]\in g\tens g\tens g$ is the leading order part of $\phi_{KZ}$ as explained in 
\cite{BegMa}. Here $r_+$ is the symmetric part of the standard classical r-matrix of $g$ and is a multiple of the Killing form. In this case the left-invariant trivector field $\tilde n$ given by the action of $n$ is some multiple of  the canonical 'Cartan tensor' that exists for any manifold $M$ acted upon by a semisimple Lie algebra. We see that  $\CC_q(M)$ is not the quantization of a usual Poisson manifold but of a 'quasi-Poisson' manifold. Such a weaker concept was proposed recently in \cite{AleKos} and we see that we have succeed in quantizing it using cotwisting. We see, moreover, that the quantum quasispace $\C_q(M)$ remains covariant but under the quantum group $\C_q(G)$ (viewed as a cotriangular coquasiHopf algebra). Let us note, however, one techincal difference from \cite{AleKos}; our
quasiPoisson manifold is associated to an action of $G$ not to a Poisson action of $G$ (these are not quite the same thing).

 Finally, we can apply all of this theory to $M=G$ i.e. to $A=\C(G)[[\hbar]]$ where $G$ acts on itself by translation, i.e. $A=H$ and the coaction is via the coproduct. This is the same idea as for our examples with finite groups, but now with $G$ a Lie group of a simple Lie algebra. 
 
 \begin{cor} The standard quantum groups $\C_q(G)$ have quasialgebra versions $\CC_q(G)$ as algebras in the category of $\C_q(G)$-comodules as a symmetric monoidal category. 
 \end{cor}
 \proof Here $\CC_q(G)=(\C(G)[[\hbar]])_F$ where we use the one-sided cotwist, in contrast to
 the Drinfeld cotwist which gives $(\C_q(G),\CR_F,\phi_{KZ})$ as explained above.  The former lives in the category of comodules of the latter. \eproof
 
 If one wants to be concrete, let $\{t^i{}_j\}$ be the matrix of generators of the classical Lie group $G$. These generate the classical $\C(G)$ with the usual 'matrix' coproduct. If we know $F$ then we know in particular the tensors
 \[ F^i{}_k{}^j{}_l=F(t^i{}_k,t^j{}_l),\quad (\Delta_2F)^i{}_l{}^j{}_m{}^k{}_n=F(t^i{}_l,t^j{}_m t^k{}_n),\quad (\Delta_1F)^i{}_l{}^j{}_m{}^k{}_n=F(t^i{}_l t^j{}_m,t^k{}_n)\]
 and so forth, where the product is in the classical $\C(G)[[\hbar]]$. Next, we denote the generators of $\CC_q(G)$ by
 $\{x^i{}_j\}$ say. They are the same as the $t^i{}_j$ but with a new product which then enjoys the deformed relations
 \[ x_1 x_2= F F^{-1}_{21} x_2 x_1.\]
We use here the standard notation in quantum group theory, where the numerical indices refer to the position in a tensor product. These relations have to be combined with the nonassociativity relations
 \[ x_1(x_2 x_3)F_{12}(\Delta_1 F)=F_{23}(\Delta_2 F)(x_1 x_2)x_3.\]
 The commutation relations reflect that the quasispace is braided-commutative with respect to the
 cotriangular structure $\CR_F$, while the nonassociativity relations reflect the associator $\phi_{KZ}$ 
 obtained from $F$.

 We also have a calculus, cocycle etc. on $\CC_q(G)$, with  $\Omega(\CC_q(G))=\Omega(\C(G)[[\hbar]])_F$. Choosing a matrix of invariant classical differential forms as generators of the classical calculus, one has similar 'F-matrix' formulae for the relations in the deformed calculus on $\CC_q(G)$, and so forth. Further details will be given elsewhere.


\begin{thebibliography}{}
\bibitem{Drinfeld} V.G. Drinfeld \emph{ Quasi-Hopf algebras}, 
\newblock Leningrad Math. J., 1:1419--1457, 1990.
\bibitem{Ma:tan} S. Majid, \emph{Tannaka-{K}rein theorem
for quasi{H}opf algebras and other results}, Contemp.
Math.,  134:219--232, 1992.
\bibitem{MaOeck} S. Majid and R. Oeckl \emph{Twisting of quantum differentials and the {P}lanck scale {H}opf algebra}
\newblock Commun. Math. Phys., 205:617--655, 1999.
\bibitem{BegMa} E. Beggs and S. Majid \emph{Semiclassical differential structures}, arXive math.QA/0306273.
\bibitem{Ma:book}
S. Majid, \emph{Foundations of quantum group theory,} Cambridge Univ. Press, Cambridge, UK, 1995
\bibitem{Con}
A. Connes, \emph{Noncommutative geometry}, Academic Press 1994
\bibitem{KMT}
J. Kustermans, G.J. Murphy, L. Tuset, \emph{Differential calculi over quantum groups and twisted cyclic cocycles,} J. Geom. Phys. 44:570--594, 2003.
\bibitem{AleKos} A. Alekseev, Y. Kosmann-Schwarzbach, and E. Meinrenken, \emph{Quasi-Poisson manifolds}, Canadian J. Math., 54(1):3Ð29, 2002.
\bibitem{Ma:euc} S. Majid \emph{{$q$}-{E}uclidean space and quantum {W}ick
rotation by twisting}, \newblock {J. Math. Phys.},
35:5025--5034, 1994. 
\bibitem{Sit}
A. Sitarz, \emph{Twist and spectral triples for isospectral deformations}, 
\newblock Lett. Math. Phys. 58: 69-79, 2001.
\bibitem{Var}
J. Varilly, \emph{Quantum symmetry groups of noncommutative spheres}, 
\newblock  Commun. Math. Phys. 221:511--523, 2001.
\bibitem{Ma:cli}
S. Majid, \emph{Noncommutative physics on Lie algebras,} $\mathbb{Z}^n_2$ \emph{lattices and Clifford algebras,} in Clifford Algebras: Application to Mathematics, Physics, and Engineering, ed. R. Ablamowicz,  Birkhauser, 2003, pp. 491-518; arXive hep-th/0302120.
\bibitem{AlbMa}
H. Albuquerque and S. Majid, \emph{Quasialgebra structure of the octonions,} Journal of algebra, 220:188--224,1999.
\bibitem{Mac}
S. Mac Lane, \emph{Categories for the working mathematician,} GTM, Vol.5, Springer-Verlag, Berlin/New York, 1974
\bibitem{Turaev}
V. Turaev, \emph{Quantum invariants of knots and 3-manifolds}, Walter de Gruter 1994
\bibitem{JoyStr} A. Joyal and R. Street, \emph{Braided monoidal categories},  Mathematics Report 86008,  Macquarie University, 1986.
\bibitem{GVF}
J.M. Gracia--Bondia, J.C. Varilly, H. Figueroa, \emph{Elements of noncommutative geometry,} Birkhauser, 2001.
\bibitem{Ma:diag} S. Majid, \emph{Diagrammatics of braided group gauge theory}, J. Knot Theor. Ramif. 
8:731--771, 1999.
\bibitem{Wor}
S.L. Woronowicz, \emph{Differential calculus on compact matrix pseudogroups (quantum groups),} Commun. Math. Phys., 122:125--170, 1989.   
\bibitem{Rieffel} M.A. Rieffel, \emph{Deformation quantization for actions of
$\mathbb{R}^d$},
\newblock {Memoirs AMS}, 106, 1993.
\end{thebibliography}
\end{document}